\documentclass[12pt,a4paper]{article}
\usepackage{slashbox,amssymb}
\usepackage{amsmath,latexsym,fullpage,color,times,tikz}
\usepackage{helvet,mathrsfs,pifont,bm}
\usepackage[pdftex]{hyperref}
\hypersetup{plainpages=True, pdfstartview=FitV,
colorlinks=true,linkcolor=blue,citecolor=blue}

\newtheorem{thm}{Theorem}

\def\dd#1{{\,\rm d}#1}
\def\ve{\varepsilon}

\def\qed{{\quad\rule{1mm}{3mm}\,}}
\def\pf{\noindent {\em Proof.}\ }
\def\mb#1{{\mathbf{#1}}}
\newcommand{\ds}{\displaystyle}

\begin{document}

\title{Multivariate records based on dominance}
\author{Hsien-Kuei Hwang and Tsung-Hsi Tsai \\
Institute of Statistical Science\\
Academia Sinica\\
Taipei 115\\
Taiwan}
\maketitle

\begin{abstract}

We consider three types of multivariate records in this paper and
derive the mean and the variance of their numbers for independent
and uniform random samples from two prototype regions: hypercubes
$[0,1]^d$ and $d$-dimensional simplex. Central limit theorems with
convergence rates are established when the variance tends to
infinity. Effective numerical procedures are also provided for
computing the variance constants to high degree of precision.

\end{abstract}

\section{Introduction}

While the one-dimensional records (or record-breakings,
left-to-right maxima, outstanding elements, etc.) of a given sample
have been the subject of research and development for more than six
decades, considerably less is known for multidimensional records.
One simple reason being that there is no total ordering for
multivariate data, implying no unique way of defining records in
higher dimensions. We study in this paper the stochastic properties
of three types of records based on the dominance relation under two
representative prototype models. In particular, central limit
theorems with convergence rates are proved for the number of
multivariate records when the variance tends to infinity, the major
difficulty being the asymptotics of the variance.

\paragraph{Dominance and maxima.}
A point $\mb{p}\in \mathbb{R}^d$ is said to \emph{dominate} another
point $\mb{q}\in \mathbb{R}^d$ if $\mb{p}- \mb{q}$ has only positive
coordinates, where the dimensionality $d\ge1$. Write $\mb{q} \prec
\mb{p}$ or $\mb{p} \succ\mb{q}$. The nondominated points in the set
$\{\mb{p}_1, \dots, \mb{p}_n\}$ are called \emph{maxima}. Maxima
represent one of the most natural and widely used partial orders for
multidimensional samples when $ d\ge2$, and have been thoroughly
investigated in the literature under many different guises and names
(such as admissibility, Pareto optimality, elites, efficiency,
skylines, \dots); see \cite{BDHT05,CHT09} and the references
therein.

\paragraph{Pareto records.}
A point $\mb{p}_k$ is defined to be a \emph{Pareto record} or a
\emph{nondominated record} of the sequence $\mb{p}_1,\ldots,
\mb{p}_n$ if
\[
    \mb{p}_k\nprec \mb{p}_i \text{ for all }1\le i<k.
\]
Such a record is referred to as a \emph{weak record} in
\cite{Gnedin07}, but we found this term less informative.

In addition to being one of the natural extensions of the classical
one-dimensional records, the Pareto records of a sequence of points
are also closely connected to maxima, the simplest connection being
the following bijection. If we consider the indices of the points as
an additional coordinate, then the Pareto records are exactly the
maxima in the extended space (the original one and the index-set) by
reversing the order of the indices. Conversely, if we sort a set of
points according to a fixed coordinate and use the ranks as the
indices, then the maxima are nothing but the Pareto records in the
induced space (with one dimension less); see \cite{Gnedin07}. See
also the recent paper \cite{CHT09} for the algorithmic aspects of
such connections.

More precisely, assume that $\mb{p}_1,\ldots ,\mb{p}_n$ are
independently and uniformly distributed (\emph{abbreviated as iud})
in a specified region $S$ and $\mb{q}_1,\ldots ,\mb{q}_n$ are iud in
the region $S\times [0,1]$. Then the distribution of the number of
Pareto records of the sequence $\mb{p}_1,\ldots, \mb{p}_n$ is equal
to the distribution of the number of maxima of the set
$\{\mb{q}_1,\ldots , \mb{q}_n\}$. This connection will be used later
in our analysis.

On the other hand, we also have, for any given regions, the
following relation between the expected number $\mathbb{E}[X_n]$ of
Pareto records and the expected number $\mathbb{E}[M_n]$ of maxima
of the same sample of points, say $\mb{p}_1,\ldots ,\mb{p}_n$,
\[
    \mathbb{E}[X_n]
    = \sum_{1\le k\le n} \frac{\mathbb{E}[M_k]}{k};
\]
see \cite{CHT09}.

\paragraph{Dominating records.}
Although the Pareto records are closely connected to maxima, their
probabilistic properties have been less well studied in the
literature. In contrast, the following definition of records has
received more attention.

A point $\mb{p}_k$ is defined to be a \emph{dominating record} of
the sequence $\mb{p}_1,\ldots ,\mb{p}_n$ if
\[
   \mb{p}_i\prec \mb{p}_k\text{ for all }1\le i<k.
\]
This is referred to as the \emph{strong record} in \cite{Gnedin07}
and the \emph{multiple maxima} in \cite{HH05}.

Let the number of dominating records falling in $A\subset S$ be
denoted by $Z_{A}$. Goldie and Resnick \cite{Goldie89} showed that
\[
    \mathbb{E}[Z_{A}]
    = \int_A \left( 1-\mu (D_{\mb{x}})\right)^{-1}
    \dd \mu (\mb{x}),
\]
where $D_{\mb{x}}=\{\mb{y}:\mb{y}\prec \mb{x}\}$. They also
calculated all the moments of $Z_{A}$ and derived several other
results such as the probability of the event $\{Z_{A}=0\}$ and the
covariance $\text{Cov} \left( Z_{A},Z_{B}\right)$.

In the special case when the $\mb{p}_i$'s are iud with a common
multivariate normal (non-degenerate) distribution, Gnedin
\cite{Gnedin98} proved that
\[
    \lambda_n
    := \mathbb{P}\{\mb{p}_n \text{ is a dominating
    record}\} \asymp n^{-\alpha }(\log n)^{(\alpha -\beta )/2}.
\]
for some $\alpha >1$ and $\beta \in \{2,3,\ldots ,d\}$. For finer
asymptotic estimates, see \cite{HH02}.

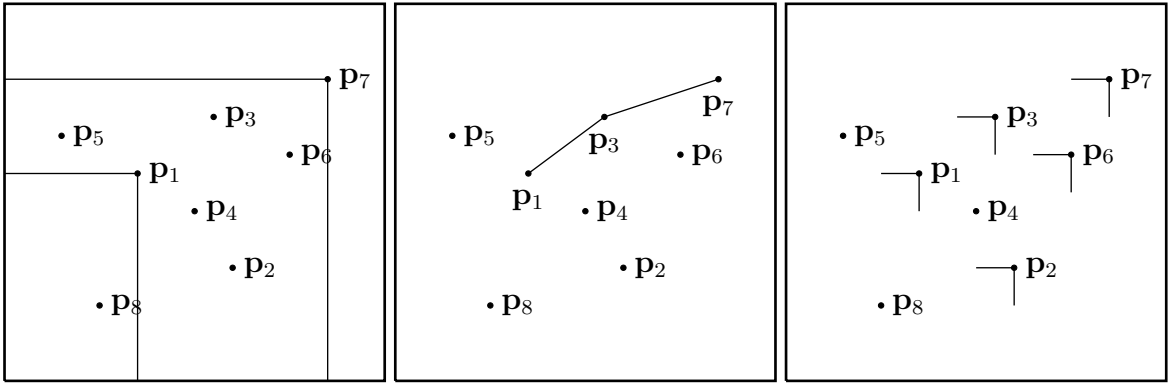
\begin{figure}[!h tbp]
\begin{center}\hspace*{0cm}\vspace*{0cm}
\begin{tikzpicture} [scale=0.5]
\draw [line width=1pt] (-7,0)--(-7,10)--(3,10)--(3,0)--(-7,0);%
\draw [fill,color=black] (-3.5,5.5) circle(2pt) node[right]{$\mb{p}_1$};%
\draw [fill,color=black] (-1,3) circle(2pt) node[right]{$\mb{p}_2$};%
\draw [fill,color=black] (-1.5,7) circle(2pt) node[right]{$\mb{p}_3$};%
\draw [fill,color=black] (-2,4.5) circle(2pt) node[right]{$\mb{p}_4$};%
\draw [fill,color=black] (-5.5,6.5) circle(2pt) node[right]{$\mb{p}_5$};%
\draw [fill,color=black] (0.5,6) circle(2pt) node[right]{$\mb{p}_6$};%
\draw [fill,color=black] (1.5,8) circle(2pt) node[right]{$\mb{p}_7$};%
\draw [fill,color=black] (-4.5,2) circle(2pt) node[right]{$\mb{p}_8$};%
\draw [line width=.5pt] (-7,5.5)--(-3.5,5.5)--(-3.5,0);%
\draw [line width=.5pt] (-7,8)--(1.5,8)--(1.5,0);%
\end{tikzpicture}\hspace*{0cm}
\begin{tikzpicture} [scale=0.5]
\draw [line width=1pt] (-7,0)--(-7,10)--(3,10)--(3,0)--(-7,0);%
\draw [fill,color=black] (-3.5,5.5) circle(2pt) node[below=2pt]{$\mb{p}_1$};%
\draw [fill,color=black] (-1,3) circle(2pt) node[right]{$\mb{p}_2$};%
\draw [fill,color=black] (-1.5,7) circle(2pt) node[below=2pt]{$\mb{p}_3$};%
\draw [fill,color=black] (-2,4.5) circle(2pt) node[right]{$\mb{p}_4$};%
\draw [fill,color=black] (-5.5,6.5) circle(2pt) node[right]{$\mb{p}_5$};%
\draw [fill,color=black] (0.5,6) circle(2pt) node[right]{$\mb{p}_6$};%
\draw [fill,color=black] (1.5,8) circle(2pt) node[below=2pt]{$\mb{p}_7$};%
\draw [fill,color=black] (-4.5,2) circle(2pt) node[right]{$\mb{p}_8$};%
\draw [line width=.5pt] (-3.5,5.5)--(-1.5,7)--(1.5,8);%
\end{tikzpicture}\hspace*{0cm}
\begin{tikzpicture} [scale=0.5]
\draw [line width=1pt] (-7,0)--(-7,10)--(3,10)--(3,0)--(-7,0);%
\draw [fill,color=black] (-3.5,5.5) circle(2pt) node[right]{$\mb{p}_1$};%
\draw [fill,color=black] (-1,3) circle(2pt) node[right]{$\mb{p}_2$};%
\draw [fill,color=black] (-1.5,7) circle(2pt) node[right]{$\mb{p}_3$};%
\draw [fill,color=black] (-2,4.5) circle(2pt) node[right]{$\mb{p}_4$};%
\draw [fill,color=black] (-5.5,6.5) circle(2pt) node[right]{$\mb{p}_5$};%
\draw [fill,color=black] (0.5,6) circle(2pt) node[right]{$\mb{p}_6$};%
\draw [fill,color=black] (1.5,8) circle(2pt) node[right]{$\mb{p}_7$};%
\draw [fill,color=black] (-4.5,2) circle(2pt) node[right]{$\mb{p}_8$};%
\draw [line width=.5pt] (-4.5,5.5)--(-3.5,5.5)--(-3.5,4.5);%
\draw [line width=.5pt] (-2,3)--(-1,3)--(-1,2);%
\draw [line width=.5pt] (-2.5,7)--(-1.5,7)--(-1.5,6);%
\draw [line width=.5pt] (-0.5,6)--(0.5,6)--(0.5,5);%
\draw [line width=.5pt] (0.5,8)--(1.5,8)--(1.5,7);%
\end{tikzpicture}
\end{center}
\caption{\emph{In this simple example, the dominating records are
$\mb{p}_1$ and $\mb{p}_7$ (left), the chain records are $\mb{p}_1$,
$\mb{p}_3$ and $\mb{p}_7$ (middle), and the Pareto records are
$\mb{p}_1$, $\mb{p}_2$, $\mb{p}_3$, $\mb{p}_6$ and $\mb{p}_7$
(right), respectively.}} \label{fg1}
\end{figure}

\paragraph{Chain records.}
Yet another type of records of multi-dimensional samples introduced
in \cite{Gnedin07} is the \emph{chain record}
\[
    \mb{p}_1\prec \mb{p}_{i_1}
    \prec \mb{p}_{i_2}\prec \cdots \prec \mb{p}_{i_k},
\]
where $1<i_1<i_2<\cdots <i_k$ and there are no $\mb{p}_{j}\mb{\
\succ p}_{i_{a}}$ with $i_{a}<j<i_{a+1}$ or $i_{a}<j\le n$. See
Figure~\ref{fg1} for an illustration of the three different types of
records.

\paragraph{Some known results and comparisons.}
If we drop the restriction of order, then the largest subset of
indices such that
\begin{equation}
    \mb{p}_{i_1}\prec \mb{p}_{i_2}
    \prec \cdots \prec \mb{p}_{i_k} \label{4}
\end{equation}
is equal to the number of maximal layers (maxima being regarded as
the first layer, the maxima of the remaining points being the
second, and so on). Assuming that $\{\mb{p}_1, \dots,\mb{p}_n\}$ are
iud in the hypercube $[0,1]^d$, Gnedin \cite{Gnedin07} proved that
the number of chain records $Y_n$ is asymptotically Gaussian with
mean and variance asymptotic to
\[
    \mathbb{E}[Y_n]
    \sim d^{-1}\log n, \quad \mathbb{V}[Y_n]\sim d^{-2}\log n;
\]
see Theorem\ref{thm-clt-cr} for an improvement.
The author also derived exact and asymptotic formulas for the
probability of a chain record $\mathbb{P}(Y_n>Y_{n-1})$ and
discussed some point-process scaling limits.

The behavior of the record sequence (\ref{4}) in $\mathbb{R}^2$ are
studied in Goldie and Resnick \cite{Goldie95}, Deuschel and Zeitouni
\cite{DZ95}. The position of the points converges in probability to
a (or a set of) deterministic curve(s). Deuschel and Zeitouni
\cite{DZ95} also proved a weak law of large number for the longest
increasing subsequence, extending a result by Vershik and Kerov
\cite{VK77} to a non-uniform setting; see also the breakthrough
paper \cite{BDJ99}. A completely different type of multivariate
records based on convex hulls was discussed in \cite{Kaluszka95}.

Chain records can in some sense be regarded as uni-directional
Parero records, and thus lacks the multi-directional feature of
Pareto records. The asymptotic analysis of the moments is in general
simpler than that for the Pareto records. On the other hand, it is
also this aspect that the chain records reflect better the
properties exhibited by the one-dimensional records. Interestingly,
the chain records correspond to the ``left-arm" (starting from the
root by always choosing the subtree corresponding to the first
quadrant) of quadtrees; see \cite{CFH07,FLLS95} and the references
therein.

\begin{figure}[!h]
\begin{center}
\begin{tikzpicture}[
 s1/.style={
 rectangle,
 rounded corners=3mm,
 minimum size=6mm,
 very thick,
 draw=white!50!black!50,
 top color=white,
 bottom color=white!50!black!20,
 scale=.9
 }]
\node[s1] (n1) at (0,0) {\large{DOMINANCE}};
\node[s1] (n2) at (-3,-1) {\large{MAXIMA}};
\node[s1] (n3) at (3,-1) {\large{RECORDS}};
\path (n1) edge[-] (n2);
\path (n1) edge[-] (n3);
\node[s1,text badly centered, text width = 2cm] (n4) at (0,-1.5)
{Pereto records};
\path (n2) edge[-] (n4);
\path (n3) edge[-] (n4);
\node[s1,text badly centered, text width = 1.8cm] (n5) at (.9,-2.8)
{chain records};
\path (n3) edge[-] (n5);
\node[s1,text badly centered, text width = 2.2cm] (n6) at (3.2,-2.8)
{dominating records};
\path (n3) edge[-] (n6);
\node (n7) at (4.8,-2.5) {\dots};
\path (n3) edge[-] (n7);
\node[s1,text badly centered, text width = 1.7cm] (n9) at (-3,-2.2)
{\large{maximal layers}};
\path (n2) edge[-] (n9);
\node[s1,text badly centered, text width = 1.5cm] (n10) at (-3,-3.5)
{\large{depth}};
\path (n9) edge[-] (n10);
\end{tikzpicture}
\end{center}
\caption{\emph{A diagram illustrating the diverse notions defined on
dominance; in particular, the Pareto records can be regarded as a
good bridge between maxima and multivariate records.}}
\end{figure}
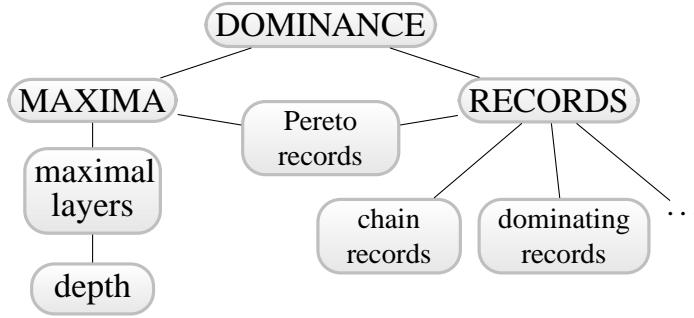

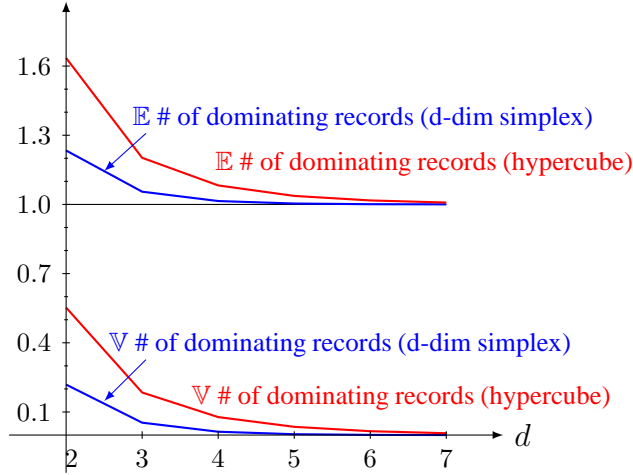
\begin{figure}[!h]
\begin{center}
\begin{tikzpicture}
\newcommand{\FONTSIZE}{\fontsize{10pt}{\baselineskip}\selectfont}
\definecolor{gen}{rgb}{0,0,0}
\draw[gen] plot coordinates{(0.000, 3.058)(5.000, 3.058)};
\definecolor{gen}{rgb}{0,0,1}
\draw[gen,thick] plot coordinates{(0.000, 3.773)(1.000,
3.227)(2.000, 3.103)(3.000, 3.070)(4.000, 3.061)(5.000, 3.059)};
\draw[gen,thick] plot coordinates{(0.000, 0.669)(1.000,
0.164)(2.000, 0.044)(3.000, 0.011)(4.000, 0.002)(5.000, 0.000)};
\definecolor{gen}{rgb}{1,0,0}
\draw[gen,thick] plot coordinates{(0.000, 5.000)(1.000,
3.676)(2.000, 3.310)(3.000, 3.171)(4.000, 3.111)(5.000, 3.083)};
\draw[gen,thick] plot coordinates{(0.000, 1.690)(1.000,
0.564)(2.000, 0.238)(3.000, 0.109)(4.000, 0.051)(5.000, 0.024)};
\draw[-latex] (-0.750 , 0.000) -- (5.750 , 0.000) node[right] {$d$};
\draw (1 , 0.05) -- (1 , -0.050) node[below] %
{\FONTSIZE{$3$}} (2, 0.05) -- (2 , -0.050)
node[below]{\FONTSIZE{$4$}} (3 , 0.05) -- (3 , -0.050)
node[below]{\FONTSIZE{$5$}} (4 , 0.05) -- (4 , -0.050)
node[below]{\FONTSIZE{$6$}} (5 , 0.05) -- (5 , -0.050)
node[below]{\FONTSIZE{$7$}}; %
\draw (0.07 , -0.050) node[below] {\FONTSIZE{$2$}};%
\draw[-latex] (0.000 , -0.50) -- (0.000 , 5.750)node[above] {};
\foreach \y/\ytext in {
0.6116,0.9174,1.529,1.8348,2.4464,2.7522,3.3638,3.6696,4.2812,4.587}
\draw[shift={(0.000, \y)}] (0.025, 0) -- (-0.025, 0);%
\draw (0.050 , 0.3058) -- (-0.050 , 0.3058)node[left]
{\FONTSIZE{$0.1$}} (0.050 , 1.2232) -- (-0.050 , 1.2232)
node[left]{\FONTSIZE{$0.4$}} (0.050 , 2.1406) -- (-0.050 , 2.1406)
node[left]{\FONTSIZE{$0.7$}} (0.050 , 3.058) -- (-0.050 , 3.058)
node[left]{\FONTSIZE{$1.0$}} (0.050 , 3.9754) -- (-0.050 , 3.9754)
node[left]{\FONTSIZE{$1.3$}} (0.050 , 4.8928) -- (-0.050 , 4.8928)
node[left]{\FONTSIZE{$1.6$}};%
\draw (4.7 , 3.6) node[red] {\FONTSIZE{$\mathbb{E}$ \#\ of
dominating records (hypercube)}};%
\draw (3.9 , 4.2) node[blue] {\FONTSIZE{$\mathbb{E}$ \#\ of
dominating records (d-dim simplex)}};%
\draw (4.4 , 0.5) node[red] {\FONTSIZE{$\mathbb{V}$ \#\ of
dominating records (hypercube)}}; %
\draw (3.6 , 1.2) node[blue] {\FONTSIZE{$\mathbb{V}$ \#\ of
dominating records (d-dim simplex)}};%
\draw[-latex][blue] (1.05 , 4) -- (0.5 , 3.5) ; \draw[-latex][blue]
(1.05 , 1) -- (0.5 , 0.45) ;
\end{tikzpicture}
\end{center}
\caption{\emph{The mean and the variance of the number of dominating
records in low dimensional random samples. In each model, the
expected number approaches $1$ very fast as $d$ increases with the
corresponding variance tending to zero.}} \label{fig-dom-rec}
\end{figure}

\paragraph{A summary of results.}
We consider in the paper the distributional aspect of the above
three types of records in two typical cases when the $\mb{p}_i$'s
are iud in the hypercube $[0,1]^d$ and in the $d$-dimensional
simplex, respectively. Briefly, hypercubes correspond to situations
when the coordinates are independent, while the $d$-dimensional
simplex to that when the coordinates are to some extent negatively
correlated. The hypercube case has already been studied in
\cite{Gnedin07}; we will discuss this briefly by a very different
approach. In addition to the asymptotic normality for the number of
Pareto records in the $d$-dimensional simplex, our main results are
summarized in the following table, where we list the asymptotics of
the mean (first entry) and the variance (second entry) in each case.
\begin{center}
\begin{tabular}{|c||c|c|}\hline
\backslashbox{Records}{Models}
& Hypercube $[0,1]^d$ & $d$-dimensional simplex \\ \hline\hline%
Dominating records & $H_n^{(d)},H_n^{(d)}-H_n^{(2d)}$ &
(\ref{mean-dom-rec}), (\ref{var-dom-rec}) \\ \hline%
Chain records &$\ds \frac1d\log n, \frac1{d^2}\log n$
\cite{Gnedin07} &
$\ds \frac1{dH_d} \log n,\frac{H_d^{(2)}}{dH_d^{3}}\log n$ \\
\hline%
Pareto records & $\ds  \frac1{d!}\left( \log n\right) ^d,\left(
 \frac1{d!} +\kappa _{d+1}\right) \left( \log n\right) ^d$
\cite{Gnedin07} &
$m_d n^{(d-1)/d}$, $v_d n^{(d-1)/d}$ \\ \hline\hline
Maxima & $=$ Pareto records in $[0,1]^{d-1}$ \cite{Gnedin07} &
$\tilde{m}_d n^{(d-1)/d}$, $\tilde{v}_dn^{(d-1)/d}$ \\ \hline
\end{tabular}
\end{center}
Here $H_{b}^{(a)}=\sum_{i=1}^{b}i^{-a}$, $\kappa_d$ is a constant
(see \cite{BDHT05}), $m_d := \frac d{d-1}\Gamma \left(
\frac1d\right)$, $v_d$ is defined in (\ref{var-Rn}), $\tilde{m}_d :=
\Gamma(\frac1d)$, $\tilde{v}_d$ is given in (\ref{cd}), and both
(\ref{mean-dom-rec}) and (\ref{var-dom-rec}) are bounded in $n$ and
in $d$; see Figure~\ref{fig-dom-rec}.

From this table, we see clearly that the three types of records
behave very differently, although they coincide when $d=1$. Roughly,
the number of dominating records is bounded (indeed less than two on
average) in both models, while the chain records have a typical
logarithmic quantity; and it is the Pareto records that reflect
better the variations of the underlying models.

\paragraph{Organization of the paper.}
We derive asymptotic approximations to the mean and the variance for
the number of Pareto records in the next section. Since the
expression for the leading coefficient of the asymptotic variance is
very messy, we then address in Section~\ref{sec-nc} the numerical
aspect of this constant. The tools we used turn out to be also
useful for several other constants of similar nature, which we
briefly discuss. We then discuss the chain records and the
dominating records.

\section{Asymptotics of the number of Pareto records}

Let
\[
    S_d
    :=\{\mb{x}:x_i\ge 0\; \text{and}\;
    0\le \left\| \mb{x}\right\| \le 1\}
\]
denote the $d$-dimensional simplex, where $\left\| \mb{x}\right\|
:=x_1+\cdots +x_d$. Assume that $\mb{p}_1,\ldots ,\mb{p}_n$ are iud
in $S_d$. Let $X_n$ denote the number of Pareto records of
$\{\mb{p}_1,\ldots , \mb{p}_n\}$. We derive in this section
asymptotic approximations to the mean and the variance and a
Berry-Esseen bound for $X_n$. The same method of proof also applies
to the number of maxima, denoted by $M_n$, which we will briefly
discuss.

Let $\mb{q}_1,\dots,\mb{q}_n$ be iud in $S_d\times [0,1]$. As
discussed in Introduction, the distribution of $X_n$ is equivalent
to the distribution of the number of maxima of
$\{\mb{q}_1,\dots,\mb{q}_n\}$.

For notational convenience, denote by $a_n\simeq b_n$ if $a_n=
b_n+O\left(n^{-1/d}\right)$.

\begin{thm} The mean and the variance of the number of Pareto
maxima $X_n$ in random samples from the $d$-dimensional simplex
satisfy \label{thm1}
\begin{align}
    \begin{split}
    \mathbb{E}[X_n]
    &\simeq n^{1-1/d}\sum_{0\le j\le d-2}\binom{d-1}{j}(-1)^{j}
    \Gamma\left( \frac{j+1}d\right)\frac d{d-1-j}\,n^{-j/d}\\
    &\qquad +(-1)^{d-1}\left(\log n+\gamma \right) ,
    \end{split} \label{mu-Rn}\\
    \mathbb{V}[X_n]
    &=\left(v_d+o(1)\right) n^{1-1/d}, \nonumber
\end{align}
where
\begin{align}
\begin{split} \label{var-Rn}
    v_d
    &:=\frac d{d-1}\Gamma \left( \frac1d\right)
    +2d^2(d-1)\sum_{1\le \ell <d}\binom d{\ell}\binom{d-2}{\ell-1}\\
    &\quad\times \int_0^1\!\!\int_v^1\!\!\int_0^\infty\!\!
    \int_0^{\infty }\!\!\int_0^\infty y^{d-\ell -1}
    w^{\ell-1}e^{-u \left( x+y\right)^d-v \left( x+w\right)^d}
    \left( e^{v x^d}-1\right) \dd w
    \dd y\dd x\dd u \dd v \\
    &\quad+2d^2\int_0^1\!\!\int_v^1\!\!\int_0^\infty \!\!\int_0^\infty
    w^{d-1}e^{-u x^d-v \left( x+w\right)^d}
    \left(e^{v x^d}-1\right)
    \dd w\dd x\dd u \dd v \\
    &\quad-2d^2\int_0^1\!\!\int_v^1\!\!\int_0^\infty\!\!
    \int_0^\infty y^{d-1}e^{-u \left( x+y\right)^d -v x^d}
    \dd y\dd x \dd u \dd v .
\end{split}
\end{align}
\end{thm}
\pf The method of proof is similar to that given in \cite{BDHT05}, 
but the technicalities are more involved. We start with the
expected value of $X_n$. Let $G_i=1_{\{\mb{q}_i \text{ is a
maxima}\}}$.
\begin{align*}
    \mathbb{E}[X_n]
    &= n\mathbb{E}[G_1]\\
    & =d!n\int_0^1\!\!\int_{S_d}
    \left( 1-z(1-\left\| \mb{x}\right\|)^d\right)^{n-1}
    \dd \mb{x}\dd z \\
    & \simeq d!n\int_0^1\!\!\int_{S_d}
    e^{-nz(1-\left\| \mb{x}\right\|)^d}
    \dd \mb{x}\dd z \\
    & =dn\int_0^1\!\!\int_0^1e^{-nz(1-y)^d}y^{d-1}
    \dd y\dd z\quad
    (y\mapsto \left\| \mb{x}\right\| ) \\
    & =dn\sum_{0\le j<d}\binom{d-1}{j}(-1)^{j}\int_0^1
    \!\!\int_0^1e^{-nzy^d}y^{j}\dd y\dd z \\
    & =\sum_{0\le j<d} n^{(d-1-j)/d}\binom{d-1}{j}(-1)^{j}
    \int_0^1\!\!\int_0^{nz}e^{-x}x^{(1+j-d)/d}z^{-(j+1)/d}
    \dd x\dd z \\
    & \simeq \sum_{0\le j\le d-2}\binom{d-1}{j}(-1)^{j}
    \Gamma \left( \frac{j+1}d\right)
    \frac d{d-1-j}\,n^{(d-1-j)/d} \\
    &\qquad +(-1)^{d-1}\left( \log n+\gamma\right) .
\end{align*}
This proves (\ref{mu-Rn}).

For the variance, we start from the second moment, which is given by
\[
    \mathbb{E}\left[ X_n^2\right]
    =\mathbb{E}[X_n]+n(n-1)\mathbb{E}\left[G_1G_2\right].
\]
Let $\mb{A}$ be the region in $\mathbb{R}^d\times[0,1]$ such that
$\mb{q}_1$ and $\mb{q}_2$ are incomparable (neither dominating the
other). Write $\mb{q}_1=(\mb{x},u)$, $\mb{q}_2=(\mb{y},v)$, $\left\|
\mb{x}\right\| _{*}:=(\left\| \mb{x}\right\| \wedge1)$ and
\[
    \mb{x}\vee \mb{y}
    :=\left( x_1\vee y_1,\cdots ,x_d\vee y_d\right) .
\]
Then by standard majorization techniques (see \cite{BDHT05})
\begin{align*}
    & n(n-1)\mathbb{E}\left[G_1G_2\right] \\
    & =n(n-1)d!^2\int_{\mb{A}}\left( 1-u (1-\left\| \mb{x}
    \right\| )^d-v (1-\left\| \mb{y}\right\| )^d
    +(u\wedge v)(1-\left\|\mb{x}\vee \mb{y}\right\|_{*})^d\right)^{n-2}
    \dd \mb{x}\dd \mb{y}\dd u \dd v \\
    & \simeq n^2d!^2\int_{\mb{A}}e^{-n[u (1-\left\| \mb{x}
    \right\| )^d+v (1-\left\| \mb{y}\right\| )^d]}
    \dd \mb{x}\dd \mb{y}\dd u \dd v \\
    & \qquad+ n^2d!^2\int_{\mb{A}}e^{-n[u (1-\left\| \mb{x}
    \right\| )^d+v (1-\left\| \mb{y}\right\| )^d]}\left(
    e^{n(u\wedge v)(1-\left\| \mb{x}\vee \mb{y}
    \right\|_{*})^d}-1\right) \dd \mb{x}
    \dd \mb{y}\dd u \dd v \\
    & \simeq \mathbb{E}\left[ X_n^2\right] \mb{-}
    J_{n,0}+\sum_{1\le \ell <d}\binom d{\ell }J_{n,\ell}+J_{n,d},
\end{align*}
where
\begin{align*}
    J_{n,0}
    & =2n^2d!^2\int_0^1\!\!\int_v^1\!\!
    \int\limits_{\substack{\mb{x}\prec \mb{y}\\ \mb{x},\mb{y}\in
    S_d}}e^{-n[u (1-\left\| \mb{x} \right\|)^d
    +v (1-\left\|\mb{y}\right\| )^d]}\dd \mb{x}
    \dd \mb{y}\dd u \dd v , \\
    J_{n,\ell}&
    =n^2d!^2\int_0^1\!\!\int_0^1\!\!
    \int\limits_{\substack{x_i>y_i, 1\le i\le \ell\\
    x_i<y_i,\ell<i\le d \\ \mb{x},\mb{y}\in S_d}}
    e^{-n[u (1-\left\| \mb{x}
    \right\| )^d+v (1-\left\|\mb{ \ y}\right\| )^d]}
    \left( e^{n(u\wedge v)(1-\left\| \mb{x}\vee
    \mb{y} \right\| _{*})^d}-1\right)
    \dd \mb{x}\dd \mb{y}\dd u \dd v, \\
    J_{n,d}& =2n^2d!^2\int_0^1\!\!\int_v^1\!\!
    \int\limits_{\substack{\mb{y}\prec\mb{x}\\
    \mb{x},\mb{y}\in S_d}}
    e^{-n[u (1-\left\| \mb{x}\right\| )^d
    +v (1-\left\| \mb{y}\right\| )^d]}\left(
    e^{n(u\wedge v)(1-\left\| \mb{x}\vee \mb{y}\right\| _{*})^d}
    -1\right) \dd \mb{x}\dd
    \mb{y}\dd u \dd z,
\end{align*}
for $1\le \ell \le d-1$.

Consider first $J_{n,\ell}$, $1\le \ell <d$. We proceed by four
changes of variables to simplify the integral starting from
\[
\begin{cases}
    x_i\mapsto \xi_i,\quad y_i\mapsto \xi_i(1-\eta_i), &
    \text{for }1\le i\le \ell; \\
    x_i\mapsto \xi_i(1-\eta_i),\quad y_i\mapsto \xi_i, &
    \text{for } \ell<i\le d ,
\end{cases}
\]
which leads to
\begin{align*}
    J_{n,\ell}
    &=(nd!)^2\int_0^1\!\!\int_0^1\!\!\int_{S_d}\!
    \int_{[0,1]^d}e^{-n\left[u
    \left( 1-\sum \xi_i+\sum^{\prime \prime }\xi_i\eta_i\right) ^d
    +v \left(1-\sum \xi_i+
    \sum^{\prime }\xi_i\eta_i\right)^d\right] } \\
    &\qquad \times \left( e^{n(u\wedge v)
    \left( 1-\sum \xi_i\right)^d}-1\right)\left(\prod \xi_i\right)
    \dd\bm{\xi}\dd \bm{\eta}\dd u\dd v ,
\end{align*}
where $\sum \xi_i := \sum_{i=1}^d\xi_i$, $\prod \xi_i :=
\prod_{i=1}^d\xi_i$, $\sum^{\prime}x_i:=\sum_{i=1}^{\ell} x_i$ and
$\sum^{\prime\prime }x_i:=\sum_{i=\ell +1}^dx_i$.

Next, by the change of variables
\[
    \xi_i\mapsto \frac1d-\xi_in^{-1/d},
    \quad \eta_i\mapsto d\eta_in^{-1/d},
\]
we have
\begin{align*}
    J_{n,\ell}
    &=d!^2\int_0^1\!\!\int_0^1\!\!\int_{S_d(n)}\!
    \int_{[0,n^{1/d}/d]^d}e^{-\left[ u
    \left( \sum \xi_i+\sum^{\prime \prime}
    \eta_i\left( 1-d\xi_i n^{-1/d}
    \right) \right) ^d+v \left( \sum \xi_i
    +\sum^{\prime }\eta_i
    \left( 1-d\xi_i n^{-1/d}\right) \right)^d\right] } \\
    &\qquad\times \left( e^{(u\wedge v)
    \left( \sum \xi_i\right) ^d}-1\right) \prod \left( 1-
    d\xi_i n^{-1/d}\right) \dd \bm{\xi}
    \dd \bm{\eta}\dd u \dd v ,
\end{align*}
where $S_d(n)=\{\bm{\xi}:\xi_i\le n^{1/d}/d$ and $\left\| \bm{\xi}
\right\| >0\}$.

We then perform the change of variables
\[
   \eta_i\mapsto \eta_i\left( 1-d\xi_i n^{-1/d}\right),
\]
and obtain
\begin{align*}
    J_{n,\ell}
    &=(d!)^2\int_0^1\int_0^1\int_{S_d(n)}
    \int_{[0,n^{1/d}/d]^d}e^{-\left[ u
    \left( \sum \xi_i+\sum^{\prime \prime}\eta_i\right)^d
    +v \left( \sum \xi_i+\sum^{\prime }\eta_i\right)^d\right] } \\
    &\qquad\times \left( e^{(u\wedge v)
    \left( \sum \xi_i\right)^d}-1\right)
    \dd \bm{\xi}\dd \bm{\eta} \dd u \dd v .
\end{align*}

Finally, we ``linearize'' the integrals by the change of variables
\[
    x\mapsto \sum \xi_i,\quad y\mapsto
    \sum\nolimits^{\prime \prime }\eta_i,\quad
    w\mapsto \sum\nolimits^{\prime }\eta_i,
\]
and get
\begin{align*}
    J_{n,\ell}
    &\simeq \frac{d\cdot d!}
    {(d-\ell -1)!(\ell -1)!}\,n^{1-1/d}
    \int_0^1\!\!\int_0^1\!\!
    \int_0^{\infty}\!\!\int_0^{\infty}\!\!
    \int_0^{\infty}y^{d-\ell -1}w^{\ell-1}
    e^{-u \left( x+y\right) ^d-v \left( x+w\right) ^d} \\
    &\qquad \times \left( e^{(u\wedge v)x^d}-1\right)
    \dd w\dd y\dd x\dd u \dd v,
\end{align*}
since the change of variables produces the factors
\[
    \frac{n^{1-1/d}}{(d-1)!},\frac{y^{d-\ell -1}}
    {(d-\ell -1)!}\text{ and }\frac{w^{\ell -1}}{(\ell -1)!}.
\]
Now by symmetry, we have
\begin{align*}
    \sum_{1\le \ell <d}\binom d{\ell }J_{n,\ell}
    &\simeq \sum_{1\le \ell <d}\binom d{\ell }
    \frac{2d\cdot d!}{(d-\ell -1)!(\ell-1)!}\,n^{1-1/d}
    \int_0^1\!\!\int_v^1\!\!\int_0^{\infty }\!\!\int_0^\infty\!\!
    \int_0^{\infty }y^{d-\ell -1}w^{\ell -1}\\
    &\qquad\qquad \times
    e^{-u \left( x+y\right)^d-v \left( x+w\right) ^d}
    \left( e^{v x^d}-1\right)
    \dd w\dd y\dd x \dd u \dd v.
\end{align*}
Proceeding in a similar manner for $J_{n,d}$, we deduce that
\begin{align*}
    J_{n,d}
    &=2d!^2\int_0^1\!\!\int_v^1\!\!\int_{S_d(n)}\!
    \int_{[0,n^{1/d}/d]^d}e^{-\left[
    u \left( \sum \xi_i\right) ^d+v
    \left(\sum \xi_i+\sum\eta_i\right) ^d\right] }
    \left( e^{(u\wedge v)\left( \sum \xi_i\right) ^d}-1\right)\!
    \dd \bm{\xi}\dd \bm{\eta}\dd u \dd v .
\end{align*}
By the change of variables $x\mapsto \sum \xi_i$,
$w\mapsto \sum\eta_i$, we have
\begin{align*}
    J_{n,d}
    &\simeq \frac{2d!^2}{((d-1)!)^2}\,n^{1-1/d}\int_0^1\!\!
    \int_v^1\!\!\int_0^{\infty }\!\!\int_0^{\infty }
    w^{d-1}e^{-u x^d-v \left(x+w\right) ^d}
    \left( e^{(u\wedge v)x^d}-1\right) \dd w\dd x\dd u \dd v \\
    &=2d^2n^{1-1/d}\int_0^1\!\!\int_v^1\!\!\int_0^{\infty }\!\!
    \int_0^\infty w^{d-1}e^{-u x^d-v
    \left( x+w\right) ^d}\left( e^{v x^d}-1\right)
    \dd w\dd x\dd u \dd v .
\end{align*}

Similarly, for $J_{n,0}$, we get
\[
    J_{n,0}
    =2d!^2\int_0^1\!\!\int_v^1\!\!\int_{S_d(n)}\!
    \int_{[0,n^{1/d}/d]^d}e^{-\left[ u \left(
    \sum \xi_i+\sum\eta_i\right) ^d
    +v \left( \sum \xi_i\right) ^d\right] }\dd
    \bm{\xi}\dd \bm{\eta}\dd u \dd v .
\]
The change of variables $x\mapsto \sum \xi_i$, $y\mapsto \sum\eta_i$
then yields
\[
    J_{n,0}
    \simeq 2d^2n^{1-1/d}\int_0^1\!\!\int_v^1\!\!
    \int_0^{\infty }\!\!\int_0^\infty y^{d-1}
    e^{-u \left( x+y\right) ^d-v x^d}
    \dd y\dd x\dd u \dd v .
\]
This completes the proof of the theorem. \qed

\paragraph{Remark.}
By the same arguments, we derive the following asymptotic estimates
for the number of maxima in $S_d$.
\begin{align*}
    \mathbb{E}[M_n]
    & \simeq \sum_{0\le j<d}\binom{d-1}{j}(-1)^{j}\Gamma
    \left( \frac{j+1}d\right) n^{(d-1-j)/d}, \\
    \mathbb{V[}M_n]
    & =\left( \tilde{v}_d+o(1)\right) n^{1-1/d},
\end{align*}
where
\begin{align} \label{cd}
\begin{split}
    \tilde{v}_d
    & =\Gamma \left(  \frac1d\right)
    +\sum_{1\le k <d}\binom dk
    \frac{dd!}{(d-k-1)!(k-1)!} \\
    & \qquad \times \int_0^{\infty }\!\!\int_0^{\infty }\!\!
    \int_0^{\infty }y^{d-k-1}w^{k-1}
    e^{-\left( x+y\right) ^d-\left( x+w\right) ^d}
    \left(e^{x^d}-1\right) \dd w
    \dd y\dd x{\,} \\
    &\qquad -2d^2\int_0^{\infty }\!\!\int_0^\infty
    y^{d-1}e^{-x^d-(x+y)^d}\dd x\dd y.
\end{split}
\end{align}

\begin{thm} The number of Pareto records in iud samples from
$d$-dimensional simplex is asymptotically normal with a rate given
by
\begin{align} \label{clt-Rn}
    \sup_x\left| \mathbb{P}\left( \frac{X_n-\mathbb{E}[X_n]}
    {\sqrt{\mathbb{V}[X_n]}}<x\right) -\Phi (x)\right|
    =O\left( n^{-(d-1)/(4d)}(\log n)^2
    +n^{-1/d}(\log n)^{1/d}\right) ,
\end{align}
where $\Phi(x)$ denotes the standard normal distribution.
\end{thm}
\pf Define the region
\[
    D_n
    :=\left\{ (\mb{x},z):\mb{x}\in S_d,z\in [0,1]\text{ and }
    z\left( 1-\left\| \mb{x}\right\| \right) ^d
    \le \frac{2\log n}n\right\}.
\]
Let $\overline{X}_n$ denote the number of maxima in $D_n$ and
$\widetilde{X}_n$ the number of maxima of a Poisson process on $D_n$
with intensity $d!n$. Then
\begin{align} \label{Rn-BE}
\begin{split}
    \left| \mathbb{P}\left(
    \frac{X_n-\mathbb{E}[X_n]}{\sqrt{\mathbb{V}[X_n]}}
    <x\right) -\Phi (x)\right|
    &\le \left| \mathbb{P}\left(
    \frac{X_n-\mathbb{E}[X_n]}{\sqrt{\mathbb{V}[X_n]} }<x\right)
    -\mathbb{P}\left(
    \frac{\overline{X}_n-\mathbb{E}[X_n]}{\sqrt{\mathbb{V
    }[X_n]}}<x\right) \right| \\
    &\quad+\left| \mathbb{P}\left(
    \frac{\overline{X}_n-
    \mathbb{E}[X_n]}{\sqrt{\mathbb{V}[X_n]}}<x\right)
    -\mathbb{P}\left( \frac{
    \widetilde{X}_n-\mathbb{E}[X_n]}
    {\sqrt{\mathbb{V}[X_n]}}<x\right) \right| \\
    &\quad+\left| \mathbb{P}\left(
    \frac{\widetilde{X}_n-\mathbb{E}[\widetilde{X}_n]}{
    \sqrt{\mathbb{V}[\widetilde{X}_n]}}<y\right) -\Phi (y)\right|
    +\left| \Phi (y)-\Phi (x)\right|,
\end{split}
\end{align}
for $x\in\mathbb{R}$, where
\[
    y=x\sqrt{\frac{\mathbb{V}[X_n]}{\mathbb{V}
    [\widetilde{X}_n]}}+\frac{\mathbb{E}
    [X_n]-\mathbb{E}[\widetilde{X}_n]}
    {\sqrt{\mathbb{V}[\widetilde{X}_n]}}.
\]
We prove that the four terms on the right-hand side of (\ref{Rn-BE})
all satisfy the $O$-bound in (\ref{clt-Rn}). For the first term, we
consider the probability
\begin{align*}
    \mathbb{P}\left( X_n\neq \overline{X}_n\right)
    &\le n\mathbb{P}\left( \mb{q}_1\notin D_n\text{ and }
    \mb{q}_1\text{is a maxima} \right) \\
    &=nd!\int_{S_d\times [0,1]-D_n}\left( 1-z(1-\left\|
    \mb{x}\right\|)^d\right) ^{n-1}
    \dd \mb{x}\dd z \\
    &\le nd!\int_{S_d\times [0,1]-D_n}
    \left( 1-\frac{2\log n}{n}\right)^{n-1}\dd \mb{x}\dd z \\
    &\le O(n^{-1}).
\end{align*}
For the second term on the right-hand side of (\ref{Rn-BE}), we use
a Poisson process approximation
\[
    \sup_{t}\left| \mathbb{P}\left( \overline{X}_n<t\right)
    -\mathbb{P}\left( \widetilde{X}_n<t\right) \right|
    \le O\left(\left| D_n\right| \right)
    =O\left( n^{-1/d}(\log n)^{1/d}\right).
\]
To bound the third term, we use Stein's method similar to the proof
for the case of hypercube given in \cite{BDHT05} and deduce that
\begin{align*}
    \sup_{y}\left| \mathbb{P}\left(\frac{\widetilde{X}_n
    -\mathbb{E}[\widetilde{X}_n]}{\sqrt{\mathbb{V}
    [\widetilde{X}_n]}}<y\right) -\Phi(y)\right|
    &=O\left(\frac{(\mathbb{E}[\widetilde{X}_n])^{1/2}Q_n}
    {(\mathbb{V}[\widetilde{X}_n])^{3/4}}\right)\\
    &=O\left( n^{-(d-1)/(4d)}\left( \log n\right) ^2\right) ,
\end{align*}
where $Q_n$ is the error term resulted from the dependence between
the cells decomposed and
\[
    Q_n
    =O\left( (\log n)^2\right) .
\]
Finally, the last term in (\ref{Rn-BE}) is bounded above as follows.
\begin{align*}
    \left| \Phi (y)-\Phi (x)\right|
    &=O\left( \frac{\left| \sqrt{\mathbb{V}[X_n]}-
    \sqrt{\mathbb{V}[\widetilde{X}_n]}\right| +\left|
    \mathbb{E}[X_n]-\mathbb{E}[ \widetilde{X}_n]\right|}
    {\sqrt{\mathbb{V}[\widetilde{X}_n]}}\right) \\
    &=O\left(n^{-(d+1)/(2d)}\right) .
\end{align*}
This proves (\ref{Rn-BE}). \qed

\paragraph{Remark.} By defining
\[
    D_n
    :=\left\{ \mb{x}:\mb{x}\in S_d\text{ and }\left( 1-\left\|
    \mb{x}\right\| \right) ^d\le \frac{2\log n}n\right\}
\]
instead and by applying the same arguments, we deduce the
Berry-Esseen bound for the number of maxima in iud samples from
$S_d$
\[
    \sup_{x}\left| \mathbb{P}\left(\frac{M_n-\mathbb{E}[M_n]}
    {\sqrt{\mathbb{V}[M_n] }}<x\right)-\Phi (x)\right|
    = O\left( n^{-(d-1)/(4d)}\log n+n^{-1/d}
    (\log n)^{1/d}\right) .
\]

\section{Numerical evaluations of the leading constants}
\label{sec-nc}

The leading constants $v_d$ (see (\ref{var-Rn})) and $\tilde{v}_d$
(see (\ref{cd})) appearing in the asymptotic approximations to the
variance of $X_n$ and to that of $M_n$ are not easily computed via
existing softwares. We discuss in this section more effective means
of computing their numerical values to high degree of precision. Our
approach is to first apply Mellin transforms (see \cite{FGD95}) and
derive series representations for the integrals by standard residue
calculations and then convert the series in terms of the generalized
hypergeometric functions
\[
    {}_{p}F_{q}(\alpha _1,\dots ,\alpha _{p};
    \beta _1,\dots ,\beta_{q};z)
    := \frac{\Gamma (\beta _1)\cdots
    \Gamma (\beta _{q})}{\Gamma(\alpha _1)\cdots
    \Gamma (\alpha _{p})}\sum_{j\ge 0}\frac{\Gamma
    (j+\alpha _1)\cdots \Gamma (j+\alpha _{p})}
    {\Gamma (j+\beta_1)\cdots \Gamma (j+\beta _{q})}
    \cdot \frac{ z^j}{j!}.
\]
The resulting linear combinations of hypergeometric functions can
then be computed easily to high degree of precision by any existing
symbolic softwares even with a mediocre laptop.

\paragraph{The leading constant $v_d$ of the asymptotic variance of
the $d$-dimensional Pareto records.} We consider the following
integrals
\begin{align*}
    C_d
    &=\sum_{1\le m<d}\binom d{m}\frac{(d-1)!}{(m-1)!(d-1-m)!} \\
    & \qquad \times \int_0^1\!\!\int_v^1\!\!\int_0^{\infty}
    \!\!\int_0^{\infty }\!\!\!\int_0^{\infty}y^{d-1-m}
    w^{m-1}e^{-u(x+y)^d-v(x+w)^d}\left(e^{vx^d}-1\right)
    \dd w\dd y\dd x\dd u\dd v \\
    & \qquad +\int_0^1\!\!\int_v^1\!\!\int_0^{\infty}
    \!\!\!\int_0^{\infty }w^{d-1}e^{-ux^d-v(x+w)^d}\left(
    e^{vx^d}-1\right) \dd w\dd x\dd u\dd v \\
    & \qquad -\int_0^1\!\!\int_v^1\!\!\int_0^{\infty}
    \!\!\!\int_0^{\infty}y^{d-1}
    e^{-u(x+y)^d-vx^d}\dd y\dd x\dd u\dd v \\
    & =:(d-1)\sum_{1\le m<d}\binom d{m}
    \binom{d-2}{m-1}I_{d,m}+I_{d,d}-I_{d,0}.
\end{align*}
Then $C_d$ is related to $v_d$ by $v_d=\frac d{d-1}\Gamma
(\frac1d)+2d^2C_d$. We start from the simplest one, $I_{d,0}$ and
use the integral representation for the exponential function
\[
    e^{-t}
    = \frac1{2\pi i}\int_{(c)}\Gamma (s)t^{-s}\dd s,
\]
where $c>0$, $\Re (t)>0$ and the integration path $\int_{(c)}$ is
the vertical line from $ c-i\infty $ to $c+i\infty $. Substituting
this representation into $I_{d,0}$, we obtain
\[
    I_{d,0}
    = \frac1{2\pi i}\int_{(c)}\Gamma
    (s)\int_0^1\!\!\int_v^1\!\!\int_0^{\infty}
    \!\!\!\int_0^{\infty}u^{-s}(x+y)^{-ds}y^{d-1}e^{-vx^d}
    \dd y\dd x\dd u\dd v\dd s.
\]
Making the change of variables $y\mapsto xy$ yields
\begin{align*}
    I_{d,0}
    &= \frac1{2\pi i}\int_{(c)}\Gamma(s)
    \int_0^1\!\!\int_v^1\!\!\int_0^{\infty}\!\!\!
    \int_0^{\infty}u^{-s}x^{d(1-s)}(1+y)^{-ds}y^{d-1}
    e^{-vx^d}\dd y\dd x\dd u\dd v\dd s \\
    & = \frac1{2\pi i}\int_{(c)}\Gamma(s)\!
    \int_0^1\!\!\int_v^1u^{-s}\!\!\left(
    \int_0^{\infty}y^{d-1}(1+y)^{-ds}\dd y\right)\!
    \left( \int_0^{\infty}x^{d(1-s)}e^{-vx^d}\!\dd x\right)\!\!
    \dd u\dd v\dd s \\
    & =\frac{d\Gamma (d-1)}{2\pi i}\int_{(c)}
    \frac{\Gamma (s)\Gamma(ds-d)\Gamma (1+ \frac1d-s)}
    {\Gamma (ds)(ds-1)}\dd s,
\end{align*}
where $1<c<1+ \frac1d$. Moving the integration path to the
right, one encounters the simple poles at $s=1+ \frac1d+j$ for
$j=0,1,\dots $. Summing over all residues of these simple poles
and proving that the remainder integral tends to zero, we get
\[
    I_{d,0}
    =\Gamma (d-1)\sum_{j\ge 0}\frac{(-1)^{j}\Gamma
    (j+1+ \frac1d)\Gamma (dj+1)}{(j+1)!\Gamma (dj+d+1)},
\]
where the terms converge at the rate $j^{-d-1+ \frac1d}$. This can
be expressed easily in terms of the generalized hypergeometric
functions.

An alternative integral representation can be derived for $I_{d,0}$
as follows.
\begin{align*}
    I_{d,0}
    &=\frac{\Gamma (d-1)}{\Gamma (d)}\sum_{j\ge 0}
    \frac{\Gamma(j+1+\frac{1 }d)}{\Gamma
    (j+2)}(-1)^{j}\int_0^1(1-x)^{d-1}x^{dj}\dd x \\
    & =\frac{\Gamma ( \frac1d)}{d-1}\int_0^1(1-x)^{d-1}
    \frac{1-(1+x^d)^{- \frac1d}}{x^d}\dd x,
\end{align*}
which can also be derived directly from the original multiple
integral representation and successive changes of variables (first $
u$, then $x$, then $v$, and finally $y$). In particular, for $d=2$,
\[
    I_{2,0}
    =\sqrt{\pi }\left( \sqrt{2}-1+\log 2 -\log(\sqrt{2}-1)\right) .
\]

Now we turn to $I_{d,d}$.
\[
    I_{d,d} = \int_0^1\!\!\int_v^1\!\!\int_0^\infty\!\!\!
    \int_0^\infty w^{d-1} e^{-ux^d-v(x+w)^d}
    \left(e^{vx^d}-1\right)\dd w\dd x\dd u\dd v.
\]
By the same arguments used above, we have
\begin{align*}
    I_{d,d}
    &=\frac{\Gamma(d)}{2d\pi i} \int_{(c)}
    \frac{\Gamma(s)\Gamma(ds-d)\Gamma(1+ \frac1d-s)}
    {\Gamma(ds)}I_{d,d}'\dd s,
\end{align*}
where $c>1$ and
\[
    I_{d,d}^{\prime}
    :=\int_0^1 v^{-s}\int_v^1
    \left( (u-v)^{s-1-\frac1d} -
    u^{s-1-\frac1d}\right) \dd u\dd v.
\]
To evaluate $I_{d,d}'$, assume first that $\frac1d<\Re(s)<1$, so
that
\begin{align*}
    I_{d,d}^{\prime}
    &= \int_0^1 v^{-s} \int_0^{1-v} u^{s-1-\frac1d}
    \dd u - \int_0^1 u^{s-1-\frac1d} \int_0^u v^{-s}\dd v \dd u \\
    &= \frac d{d-1} \left( \frac{\Gamma(1-s)\Gamma(s-\frac1d)}
    {\Gamma(1-\frac1d)} - \frac1{1-s}\right).
\end{align*}
Now the right-hand side is well-defined for $\frac1d<\Re(s)<2$.
Substituting this into $I_{d,d}$, we obtain
\[
    I_{d,d}
    = \frac{\Gamma(d-1)}{2\pi i} \int_{(c)}
    \frac{\Gamma(s)\Gamma(ds-d) \Gamma(1+\frac1d-s)}
    {\Gamma(ds)}\left(\frac{\Gamma(1-s)\Gamma(s-\frac1d)}
    { \Gamma(1-\frac1d)} - \frac1{1-s}\right)\dd s,
\]
where $1<c<1+\frac1d$. For computational purpose, we use the
functional equation for Gamma function
\[
    \Gamma(1-s)\Gamma(s) = \frac{\pi}{\sin\pi s},
\]
so that
\[
    I_{d,d}
    = \frac{\Gamma(d-1)}{2\pi i} \int_{(c)}
    \frac{\pi\Gamma(ds-d)}
    {\Gamma(ds)\sin(\pi(s-\frac1d))}\left(
    \frac{\pi}{\Gamma(1-\frac1d)\sin(\pi s )} +
    \frac{\Gamma(s-1)}{\Gamma(ds)
    \Gamma(s-\frac1d)}\right)\dd s.
\]
In this case, we have simples poles at $s=j+1/d$ for both integrands
and $s=j$ for the first integrand to the right of $\Re(s)=1$ for
$j=2,3,\dots$. Thus summing over all the residues and proving that
the remainder integral goes to zero, we obtain
\begin{align*}
    I_{d,d}
    =\Gamma(d-1)\Gamma(\tfrac1d) \sum_{j\ge2}
    \frac{\Gamma(dj-d)}{j\Gamma(dj)} -\Gamma(d-1)
    \sum_{j\ge2}\frac{(-1)^j\Gamma(j-1+\frac1d)
    \Gamma(dj-d+1)} {\Gamma(j)\Gamma(dj+1)}.
\end{align*}
A similar argument as that used for $I_{d,0}$ gives the alternative
integral representation
\[
    I_{d,d}
    = \frac{\Gamma(\frac1d)}{d(d-1)}
    \left(-1+ \int_0^1\left(1-t^{\frac1d}\right)^{d-1}
    \left(t^{\frac1d-1}(1+t)^{-\frac1d}+ \frac{
    -\log(1-t)-t}{t^2}\right) \dd t\right).
\]
In particular, for $d=2$,
\[
    I_{2,2}
    = \sqrt{\pi}\left(2-\sqrt{2}-2\log 2+\log(\sqrt2+1)\right).
\]

Now we consider $I_{d,m}$ for $1\le m<d$.
\[
    I_{d,m} := \int_0^1\!\!\int_v^1\!\!\int_0^\infty
    \!\!\int_0^\infty\!\!\!\int_0^\infty y^{d-1-m}w^{m-1}
    e^{-u(x+y)^d-v(x+w)^d} \left(e^{vx^d}-1\right)
    \dd w\dd y\dd x\dd u\dd v,
\]
which by the same arguments leads to
\begin{align*}
    I_{d,m}
    &= \frac1{2\pi i} \int_{(c)} \Gamma(s)
    \int_0^1\!\!\int_v^1 u^{-s}\!\!
    \left(\int_0^\infty y^{d-1-m} (1+y)^{-ds}
    \dd y\right) \\ &\qquad \times \int_0^\infty w^{m-1}
    \left(\int_0^\infty x^{d(1-s)}
    \left(e^{-vx^d((1+w)^d-1)}-e^{-vx^d(1+w)^d} \right)
    \dd x\right)\dd w\dd u\dd v \dd s \\
    &= \frac{d\Gamma(d-m)}{2(d-1)\pi i} \int_{(c)}
    \frac{\Gamma(s)\Gamma(ds-d+m) \Gamma(1+\frac1d-s)}
    {\Gamma(ds)(ds-1)}W_m(s) \dd s,
\end{align*}
where $1<c<1+\frac1d$ and
\begin{align*}
    W_m(s)
    &:= \int_0^\infty w^{m-1}
    \left(\left((1+w)^d-1\right)^{s-1-\frac1d}
    -(1+w)^{ds-d-1} \right) \dd w \\
    &= \frac1d\int_0^1 t^{-s} (t^{-\frac1d}-1)^{m-1}
    \left((1-t)^{s-1-\frac1d}-1\right)\dd t \\
    &=\frac1d\sum_{0\le \ell<m}\binom{m-1}{\ell}
    (-1)^{m-1-\ell}\left(\frac{ \pi\Gamma(s-\frac1d)}
    {\Gamma(1-\frac{\ell+1}d)\Gamma(s+\frac{\ell}d)
    \sin(\pi(s+\frac{\ell}d))}-\frac1{1-\frac\ell d-s}\right),
\end{align*}
for $\frac1d<\Re(s)<2-(m-1)/d$. Note that each term has no pole at
$s=1- \frac{\ell}d$. Thus
\[
 I_{d,m}=\frac{\Gamma(d-m)}{d-1}\sum_{0\le \ell<m}
 \binom{m-1}{\ell}(-1)^{m-1-\ell} I_{d,m,\ell},
\]
where
\begin{align*}
    I_{d,m,\ell}
    &:= \frac1{2\pi i} \int_{(c)}
    \frac{\Gamma(s)\Gamma(ds-d+m)
    \Gamma(1+\frac1d-s)} {\Gamma(ds)(ds-1)} \\
    &\qquad \times \left(\frac{\pi\Gamma(s-\frac1d)}
    {\Gamma(1-\frac{\ell+1}d )\Gamma(s+\frac{\ell}d)
    \sin(\pi(s+\frac{\ell}d))} -\frac1{1-\frac\ell d-s}\right)
    \dd s.
\end{align*}
We then deduce that the integral equals the sum of the residues at
$s=j+\frac1d$ and $s=j+1-\frac\ell d$
\begin{align*}
    I_{d,m,\ell}
    &= -\frac{\Gamma(\frac{\ell+1}d)}d \sum_{j\ge1}
    \frac{\Gamma(j+1+\frac1d)
    \Gamma(dj+m+1)} {(j+1)\Gamma(dj+d+1)
    \Gamma(j+1+\frac{\ell+1}d)} \\
    &\qquad + \frac1d\sum_{j\ge1}
    \frac{(-1)^j\Gamma(j+1+\frac1d)\Gamma(dj+m+1)
    } {(j+1)!\Gamma(dj+d+1)(j+\frac{\ell+1}d)} \\
    &\qquad + \Gamma(\tfrac{\ell+1}d)\sum_{j\ge1}
    \frac{\Gamma(j+1-\frac\ell
    d)\Gamma(dj+m-\ell)} {j!\Gamma(dj+d-\ell)(dj+d-\ell-1)} \\
    &=I_{d,m,\ell}^{[1]} + I_{d,m,\ell}^{[2]} +I_{d,m,\ell}^{[3]}.
\end{align*}
It follows that
\begin{align*}
    &C_d-I_{d,d}+I_{d,0}\\
    &\qquad= (d-1)\sum_{1\le m<d} \binom d{m}\binom{d-2}{m-1}
    I_{d,m} \\
    &\qquad= \sum_{1\le m<d} \binom d{m}\frac{(d-2)!}{(m-1)!}
    \sum_{0\le \ell<m} \binom{m-1}{\ell} (-1)^{m-1-\ell} \left(
    I_{d,m,\ell}^{[1]}+I_{d,m,\ell}^{[2]}
    + I_{d,m,\ell}^{[3]}\right) \\
    &\qquad=: C_d^{[1]} + C_d^{[2]} + C_d^{[3]}.
\end{align*}
For further simplification of these sums, we begin with $C_d^{[2]}$.
Note first that
\begin{align*}
    &\sum_{0\le \ell<m} \binom{m-1}{\ell}(-1)^{m-1-\ell}
    I_{d,m,\ell}^{[2]} \\ &\qquad = \frac1d\sum_{j\ge1}
    \frac{(-1)^j\Gamma(j+1+\frac1d)\Gamma(dj+m+1) }
    {(j+1)!\Gamma(dj+d+1)} \sum_{0\le \ell<m} \binom{m-1}{\ell}
    (-1)^{m-1-\ell}\frac1{j+\frac{\ell+1}d} \\
    &\qquad= (-1)^{m-1}(m-1)!\sum_{j\ge1}
    \frac{(-1)^j\Gamma(j+1+\frac1d) \Gamma(dj+1)}
    {(j+1)!\Gamma(dj+d+1)}.
\end{align*}
Thus
\begin{align*}
    C_d^{[2]}
    &= \sum_{1\le m<d} \binom d{m}\frac{(d-2)!}{(m-1)!}
    \sum_{0\le \ell<m}
    \binom{m-1}{\ell}(-1)^{m-1-\ell} I_{d,m,\ell}^{[2]} \\
    &=(d-2)!\sum_{1\le m<d} \binom d{m}(-1)^{m-1} \sum_{j\ge1}
    \frac{(-1)^j\Gamma(j+1+\frac1d)\Gamma(dj+1)}
    {(j+1)!\Gamma(dj+d+1)} \\
    &= (1+(-1)^d) \left(I_{d,0} -
    \frac{\Gamma(1+\frac1d)}{d(d-1)}\right).
\end{align*}
Accordingly, $C_d^{[2]}=0$ for odd values of $d$.

For the other two sums containing $I_{d,m,\ell}^{[1]}$ and
$I_{d,m,\ell}^{[3]}$, we use the identity
\[
    \sum_{\ell<m<d} \frac{(N+m)!(-1)^m} {m!(d-m)!(m-1-\ell)!}
    =\frac{(-1)^dN!} { (d-1-\ell)!}\left(\binom{N+1+\ell}d
    -\binom{N+d}d\right).
\]
Then
\begin{align*}
    C_d^{[1]}
    &=\sum_{1\le m<d} \binom d{m}\frac{(d-2)!}{(m-1)!}
    \sum_{0\le\ell<m} \binom{m-1}{\ell} (-1)^{m-1-\ell}
    I_{d,m,\ell}^{[1]} \\
    &= (d-2)!\sum_{0\le \ell\le d-2}\frac{d!}{\ell!}(-1)^\ell
    \frac{\Gamma(\frac{ \ell+1}d)}d \sum_{j\ge1}
    \frac{\Gamma(j+1+\frac1d)} {(j+1)\Gamma(dj+d+1)
    \Gamma(j+1+\frac{\ell+1}d)} \\
    &\qquad \times \sum_{\ell<m<d} \frac{\Gamma(dj+m+1)(-1)^m}
    {m!(d-m)!(m-1-\ell)!} \\
    &= \frac{(-1)^d}{d(d-1)}\sum_{0\le \ell \le d-2}
    \binom{d-1}{\ell}(-1)^\ell \Gamma(\tfrac{\ell+1}d)
    \sum_{j\ge1}\frac{\Gamma(j+1+\frac1d)}
    { (j+1)\Gamma(j+1+\frac{\ell+1}d)}\left(
    \frac{\binom{dj+\ell+1}d}{\binom{ dj+d}d} -1\right).
\end{align*}
Note that
\[
    \frac{\binom{dj+\ell+1}d}{\binom{dj+d}d}-1
    = O(j^{-1})\qquad(0\le \ell\le d-2),
\]
for large $j$, so that the series is absolutely convergent.

Similarly,
\begin{align*}
    C_d^{[3]}
    &= \sum_{1\le m<d} \binom d{m}\frac{(d-2)!}{(m-1)!}
    \sum_{0\le\ell<m} \binom{m-1}{\ell}
    (-1)^{m-1-\ell} I_{d,m,\ell}^{[3]} \\
    &= \frac{(-1)^d}{d-1}\sum_{0\le \ell \le d-2}
    \binom{d-1}{\ell}(-1)^{\ell-1} \Gamma(\tfrac{\ell+1}d)
    \sum_{j\ge1}\frac{\Gamma(j+1-\frac \ell d)} {j!(dj+d-\ell-1)}
    \left(\frac{\binom{dj}d}{\binom{dj+d-\ell-1}d} -1\right).
\end{align*}

Since $v_d=\frac d{d-1}\Gamma (\frac1d)+2d^2C_d$, we obtain, by
converting the series representations into hypergeometric functions,
the following approximate numerical values of $v_d$.
\begin{align*}
    v_2&\approx2.86126\,35493\,11178\,82531\,14379,\\
    v_3&\approx3.22524\,36444\,05576\,89660\,59392,\\
    v_4&\approx3.97797\,27442\,19455\,29292\,64760,\\
    v_5&\approx4.84527\,39171\,62611\,42226\,50057,\\
    v_6&\approx5.76349\,95321\,96568\,64812\,77416,\\
    v_7&\approx6.70865\,12250\,86590\,36364\,34742,\\
    v_8&\approx7.66955\,04435\,24665\,04704\,24808,\\
    v_9&\approx8.64032\,79742\,08287\,24931\,00067,\\
    v_{10}&\approx9.61764\,75521\,13755\,73944\,20940,\\
    v_{11}&\approx10.59949\,78766\,56951\,63098\,76869,\\
    v_{12}&\approx11.58460\,78314\,60409\,77794\,37163.
\end{align*}
In particular, $v_2$ has a closed-form expression
\[
    v_2 = \frac23\sqrt{\pi}\left(2\pi^2-9-12\log 2\right).
\]

\paragraph{The leading constant $\tilde{v}_d$ of the asymptotic
variance of the $d$-dimensional maxima.} Let
\[
    J_{d,0}:=2d^2\int_0^{\infty }\!\!\!\int_0^{\infty}
    y^{d-1}e^{-x^d-(x+y)^d}\dd x\dd y,
\]
and
\[
    J_{d,k}:=\frac{dd!}{(d-k -1)!(k-1)!}\int_0^{\infty}\!\!
    \int_0^{\infty}\!\!\int_0^{\infty}y^{d-k -1}
    w^{k-1}e^{-\left(x+y\right)^d-\left( x+w\right)^d}
    \left( e^{x^d}-1\right)\dd w\dd y\dd x.
\]
Then (see (\ref{cd}))
\[
    \tilde{v}_d
    = \Gamma\left(\frac 1d\right) +
    \sum_{1\le k<d} \binom{d}{k} J_{d,k} - J_{d,0}.
\]
Consider first $J_{d,0}$. By expanding $(1+x^d)^{-1-\frac1d}$,
interchanging and evaluating the integrals, we obtain
\begin{align*}
    J_{d,0}
    &=2\Gamma \left( {\frac1d}\right)
    \int_0^1\frac{(1-x)^{d-1}}{(1+x^d)^{1+\frac1d}}\dd x \\
    &=2d!\sum_{j\ge 0}\frac{\Gamma (j+1+\frac1d)
    \Gamma (dj+1)}{\Gamma(j+1)\Gamma (dj+d+1)}\,(-1)^j,
\end{align*}
the general terms converging at the rate $O(j^{-d-\frac1d})$. The
convergence rate can be accelerated as follows.
\begin{align*}
    J_{d,0}
    &=2\Gamma \left( 1+{\frac1d}\right)
    \int_0^1x^{\frac1d-1}(1-x^{\frac1d})^{d-1}
    (1+x)^{-1-\frac1d}\dd x \\
    & =2\Gamma \left( 1+{\frac1d}\right)
    \sum_{r\ge0}2^{-r-1-\frac1d }
    \int_0^1(1-x)^r x^{\frac1d-1}
    (1-x^{\frac1d})^{d-1}\dd x \\
    & =\Gamma \left(1+{\frac1d}\right) 2^{-\frac1d}
    \sum_{0\le \ell<d}\binom{d-1}{\ell}(-1)^\ell\Gamma
    \left( {\frac{\ell+1}d}\right)\sum_{j\ge 0}
    \frac{ \Gamma (j+1+\frac1d)}
    {\Gamma (j+1+\frac{\ell+1}d)}\,2^{-j},
\end{align*}
the convergence rate being now exponential. In terms of the
generalized hypergeometric functions, we have
\[
    J_{d,0}=\Gamma \left( {\frac1d}\right) 2^{-\frac1d}
    \sum_{0\le \ell<d}\binom{d-1 }{\ell}
    \frac{(-1)^\ell}{\ell+1}{}_2F_1\left(1+{\frac1d},
    1;1+{\ \frac{\ell+1}d};{ \frac12}\right).
\]

The integrals $J_{d,k}$ can be simplified as follows.
\begin{align*}
    J_{d,k+1}
    &=d^2(d-1)\binom{d-2}k\int_0^{\infty}
    (e^{x^d}-1)\int_{x}^{\infty }e^{-y^d}\\
    &\qquad \times\int_{x}^{\infty}(y-x)^{d-2-k}(z-x)^ke^{-z^d}
    \dd z\dd y\dd x \\
    & =2d^2(d-1)\binom{d-2}k\int_0^{\infty}
    e^{-y^d}\int_0^{y}e^{-z^d} \\
    &\qquad \times
    \int_0^{z}(e^{x^d}-1)(y-x)^{d-2-k}(z-x)^k\dd x\dd z\dd y \\
    & =2(d-1)\Gamma \left( {\frac1d}\right) \binom{d-2}k
    \int_0^1(1-x)^k\int_0^1(1-xz)^{d-2-k}z^{k+1} \\
    & \qquad \times \left(
    \frac1{(1+z^d-x^dz^d)^{1+\frac1d}}
    -\frac1{(1+z^d)^{1+\frac1d}}\right)\dd z\dd x \\
    & =J_{d,k+1}'+J_{d,k+1}''.
\end{align*}

By the same proof used for $J_{d,0}$, we have
\begin{align*}
    J_{d,k+1}^{\prime \prime }
    & =-2(d-1)\Gamma \left(
    {\frac1d}\right)\binom{ d-2}k\int_0^1(1-x)^k\\
    &\qquad \qquad\times
    \int_0^1(1-xz)^{d-2-k}z^{k+1}(1+z^d)^{-1-\frac1d}\dd z\dd x \\
    & =(-1)^{k+1}2^{-\frac1d}\Gamma \left( {\frac1d}\right)
    \sum_{k<j<d}\binom{ d-1}{j}\frac{(-1)^{j}}{j+1}\,{}_2F_1\left(
    1+{\frac1d},1;1+{\ \frac{j+1}d };{\frac12}\right) .
\end{align*}

Similarly,
\begin{align*}
    J_{d,k+1}^{\prime }
    &=2(d-1)\Gamma \left( {\frac1d}\right)
    \binom{d-2}k \int_0^1(1-x)^k\\
    &\qquad \times
    \int_0^1(1-xz)^{d-2-k}z^{k+1}
    (1+z^d-x^dz^d)^{-1-\frac1d} \dd z\dd x \\
    & =2\Gamma \left( {\frac1d}\right) (d-1)!\sum_{0\le j\le
    d-2-k}\frac{(-1)^{j}}{j!(d-2-k-j)!} \\
    & \qquad \times
    \sum_{0\le \ell \le k}\frac{(-1)^{\ell
    } }{\ell !(k-\ell )!}\cdot
    \frac{{}_3F_2(1+\frac1d,\frac{k+j+2}d,1;1+
    \frac{\ell +j+1}d,1+\frac{k+j+2}d;-1)}{(\ell +j+1)(k+j+2)}.
\end{align*}

Thus we obtain the following numerical values for the limiting
constant $\tilde{v}_d$ of $\mathbb{V[}M_n]/n^{(d-1)/d}$
\begin{align*}
    \tilde{v}_2&\approx0.68468\,89279\,50036\,17418\,09957 , \\
    \tilde{v}_3&\approx1.48217\,31873\,40583\,68601\,11369 , \\
    \tilde{v}_4&\approx2.35824\,37612\,02486\,93742\,28054 , \\
    \tilde{v}_5&\approx3.27773\,90059\,79491\,26684\,80858 , \\
    \tilde{v}_6&\approx4.22231\,09450\,77067\,79998\,34338, \\
    \tilde{v}_7&\approx5.18220\,76686\,16078\,48517\,29967, \\
    \tilde{v}_8&\approx6.15196\,29023\,77474\,45508\,28039, \\
    \tilde{v}_9&\approx7.12835\,13658\,43360\,52793\,29089, \\
    \tilde{v}_{10}&\approx8.10938\,23221\,15849\,82527\,77117,\\
    \tilde{v}_{11}&\approx9.09377\,74697\,86680\,89694\,70616, \\
    \tilde{v}_{12}&\approx10.0806\,86465\,19733\,08113\,16376.
\end{align*}
In particular, $\tilde{v}_2 = \sqrt{\pi}(2\log 2-1)$; see
\cite{BHLT01}.

\paragraph{Yet another constant in \cite{CQ97}.}
A similar but simpler integral to (\ref{cd}) appeared in \cite{CQ97},
which is of the form
\[
    K_d
    :=\int_0^{\infty }\!\!\!\int_0^\infty\!\!\!
    \int_0^\infty(u+w)^{d-2}e^{-(u+x)^d+x^d-(w+x)^d}
    \dd x\dd u\dd w,
\]
(this $K_d$ is indeed their $K_{d-1}$). By Mellin inversion formula
for $e^{-t}$, we obtain
\begin{align*}
    K_d
    & = \frac1{2\pi i}\int_{(c)}\Gamma (s)
    \int_0^\infty\!\!\!\int_0^{\infty }\!\!\!
    \int_0^\infty(u+w)^{d-2}(u+x)^{-ds}
    e^{-(w+x)^d+x^d}\dd x\dd
    u\dd w\dd s \\
    & = \frac1{2d\pi i}\int_{(c)}\Gamma (s)
    \Gamma (1+\tfrac{1}d-s)\\
    &\qquad \times \int_0^{\infty }\!\!\!
    \int_0^\infty(u+w)^{d-2}(1+u)^{-ds}
    \left( (1+w)^d-1\right)^{s-1- \frac1d}\dd u\dd w\dd s.
\end{align*}
Expanding the factor $(u+w)^{d-2}$, we obtain $K_d=\sum_{0\le m\le
d-2} \binom{d-2}{m}K_{d,m}$, where
\begin{align}
    K_{d,m}
    & := \frac1{2d\pi i}\int_{(c)}\Gamma (s)
    \Gamma (1+\tfrac{1}d-s)\left( \int_0^{\infty }
    u^{m}(1+u)^{-ds}\dd u\right) \nonumber\\
    & \qquad \times \left( \int_0^{\infty }w^{d-2-m}
    \left( (1+w)^d-1\right)^{s-1- \frac1d}
    \right) \dd w \nonumber \\
    & = \frac1{2d\pi i}\int_{(c)}\Gamma (s)
    \Gamma (1+\tfrac{1}d-s)
    B(m+1,ds-m-1)U_{m}(s)\dd s. \label{Kdm}
\end{align}
Here
\begin{align*}
    U_{m}(s)
    & :=\int_0^{\infty }w^{d-2-m}
    \left( (1+w)^d-1\right) ^{s-1- \frac1d}\dd w \\
    & = \frac1d\int_0^1t^{-s}(1-t)^{s-1- \frac1d}
    \left(t^{- \frac1d}-1\right) ^{d-2-m}\dd t \\
    & = \frac1d\sum_{0\le \ell \le d-2-m}
    \binom{d-2-m}{\ell }(-1)^{d-2-m-\ell}
    B(1-s-\tfrac{\ell }d,s-\tfrac{1}d).
\end{align*}
Thus we obtain
\begin{align*}
    K_d
    & = \frac1{d^2}\sum_{0\le m\le d-2}
    \binom{d-2}{m}\sum_{0\le \ell \le d-2-m}
    \binom{d-2-m}{\ell }(-1)^{d-2-m-\ell }m!
    \Gamma (\tfrac{\ell +1}d) \\
    & \qquad \times \sum_{j\ge 0}
    \left( \frac{\Gamma (j+1-\frac{\ell }d)
    \Gamma(dj+d-\ell -m-1)}{j!\Gamma (dj+d-\ell )}
    -\frac{\Gamma (j+1+ \frac1d)\Gamma (dj+d-m)}
    {\Gamma (j+1+\frac{\ell +1}d)\Gamma (dj+d+1)}\right) .
\end{align*}
This readily gives, by converting the above series into
hypergeometric functions, the numerical values of the first few
$K_d$,
\begin{align*}
    K_2 &\approx 0.30714\,28473\,56944\,02518\,48954,\\
    K_3 &\approx 0.21288\,24684\,73220\,99693\,80676,\\
    K_4 &\approx 0.19494\,67028\,23033\,18190\,40460,\\
    K_5 &\approx 0.20723\,21512\,99671\,45854\,93769,\\
    K_6 &\approx 0.24331\,17024\,51836\,72554\,88428,\\
    K_7 &\approx 0.30744\,56566\,07893\,22242\,37300,\\
    K_8 &\approx 0.41127\,01058\,90385\,83873\,59349,\\
    K_9 &\approx 0.57571\,68456\,67243\,64328\,08087,\\
    K_{10} &\approx 0.83615\,82236\,77116\,00233\,16115,\\
    K_{11} &\approx 1.25179\,63251\,14070\,86480\,31485,\\
    K_{12} &\approx 1.92201\,04035\,18847\,36012\,85304.
\end{align*}
These are consistent with those given in Chiu and Quine (1997). In
particular, $K_2= \frac1{4}\sqrt{\pi }\log 2$. Further
simplification of this formula can be obtained as above, but the
resulting integral expression is not much simpler than
\[
    \frac{\Gamma ( \frac1d)}{d^{4}}
    \int_0^1\!\!\int_0^1\left(u^{-  \frac1d}
    +v^{- \frac1d}-2\right)^{d-2}u^{-1- \frac1d}
    v^{-1-\frac{1 }d}\left(u^{-1}+v^{-1}-1
    \right) ^{-1- \frac1d}\dd u\dd v.
\]

\section{Asymptotics of the number of chain records}

We consider in this section the number of chain records of random
samples from $d$-dimensional simplex; the tools we use are different
from \cite{Gnedin07} and apply also to chain records for hypercube
random samples, which will be briefly discussed. For other types of
results, see \cite{Gnedin07}.

\subsection{Chain records of random samples from $d$-dimensional
simplex}

Assume that $\mb{p}_1,\ldots ,\mb{p}_n$ are iud in the
$d$-dimensional simplex $S_d$. Let $Y_n$ denote the number of chain
records of this sample. Then $Y_n$ satisfies the recurrence
\begin{align} \label{Xn-d-simplex}
    Y_ n\stackrel{d}{=} 1+Y_{I_n}\qquad(n\ge1),
\end{align}
with $Y_0:=0$, where
\[
    \mathbb{P}(I_n=k)
    = \pi_{n,k}
    = d \binom{n-1}{k} \int_0^1 t^{kd}
    (1-t^d)^{n-1-k}(1-t)^{d-1}\dd t,
\]
for $0\le k<n$. An alternative expression for the probability
distribution $\pi_{n,k}$ is
\[
    \pi_{n,k}
    = \binom{n-1}{k} \sum_{0\le j<d} \binom{d-1}{j}
    (-1)^j \frac{\Gamma(n-k)\Gamma\left(k+\frac{j+1}{d}\right)}
    {\Gamma\left(n+\frac{j+1}d\right)},
\]
which is more useful from a computational point of view.

Let
\[
    (z+1)\cdots(z+d)-d!
    = z\prod_{1\le \ell<d}(z-\lambda_\ell),
\]
where the $\lambda_\ell$'s are all complex ($\not\in\mathbb{R}$),
except when $d$ is even (in that case, $-d-1$ is the unique real
zero among $\{\lambda_1,\dots,\lambda_{d-1}\}$). Interestingly,
an essentially the same equation also arises in the analysis of
random increasing $k$-trees; see \cite{DHBS09}.

\begin{thm} The number of chain records $Y_n$ for random samples from
$d$-dimensional simplex is asymptotically normally distributed in
the following sense
\begin{align} \label{Yn-BEB}
    \sup_{x\in\mathbb{R}} \left|\mathbb{P}
    \left(\frac{Y_n-\mu_s\log n}{\sigma_s \sqrt{\log n}}<x\right)
    -\Phi(x)\right| = O\left((\log n)^{-1/2}\right),
\end{align}
where $\mu_s := 1/(dH_d)$ and $\sigma_s := \sqrt{H_d^{(2)}
/(dH_d^3)}$. The mean and the variance are asymptotic to
\begin{align}
    \mathbb{E}[Y_n]
    &= \frac{H_n}{dH_d} + c_1 + O(n^{-\ve}),\label{mu-Yn}\\
    \mathbb{V}[Y_n]
    &= \frac{H_d^{(2)}}{dH_d^3}\,H_n +c_2 +O(n^{-\ve}),
    \label{var-Yn}
\end{align}
respectively, where
\begin{align*}
    c_1
    &= \frac1{dH_d}\sum_{1\le \ell <d}
    \left(\psi\left(-\frac{\lambda_\ell}d\right)
    -\psi\left(\frac\ell d\right)\right), \\
    c_2
    &= \frac16 +\frac{\pi^2}{6d^2H_d^2}
    - \frac{2H_d^{(3)}}{3H_d^3}
    +\frac{(H_d^{(2)})^2}{2H_d^4}+ \frac1{d^2H_d^2}
    \sum_{1\le \ell <d}\left(\psi'\left(
    -\frac{\lambda_\ell}d\right)-\psi'\left(\frac\ell d
    \right)\right) \\
    &\qquad + \frac{c_1H_d^{(2)}}{H_d^2} -\frac{2d!}{H_d}
    \sum_{j\ge1} \frac{(dj+1)\cdots(dj+d)(H_{dj+d}-H_{dj})}
    {((dj+1)\cdots(dj+d)-d!)^2}.
\end{align*}
\end{thm}
The error terms in (\ref{mu-Yn}) and (\ref{var-Yn}) can be further
refined, but we content ourselves with the current forms for
simplicity.

\paragraph{Expected number of chain records.}
We begin with the proof of (\ref{mu-Yn}). Consider the mean $\mu_n
:= \mathbb{E}[Y_n]$. Then $\mu_0=0$ and, by (\ref{Xn-d-simplex}),
\begin{align}\label{mu-rr}
    \mu_n = 1+\sum_{0\le k<n} \pi_{n,k} \mu_k\qquad(n\ge1).
\end{align}
Let $\tilde{f}(z) := e^{-z} \sum_{n\ge0} \mu_n z^n/n!$ denote the
Poisson generating function of $\mu_n$. Then, by (\ref{mu-rr}),
\[
    \tilde{f}(z) + \tilde{f}'(z) = 1+ d
    \int_0^1 \tilde{f}(t^dz) (1-t)^{d-1} \dd t.
\]
Let $\tilde{f}(z) = \sum_{n\ge0}\tilde{\mu}_n z^n/n!$. Taking the
coefficients of $z^n$ on both sides gives the recurrence
\[
    \tilde{\mu}_n + \tilde{\mu}_{n+1}
    = \frac{d!}{(dn+1)\cdots(dn+d)}\,\tilde{\mu}_n
    \qquad(n\ge1).
\]
Solving this recurrence using $\tilde{\mu}_1=1$ yields
\[
    \tilde{\mu}_n
    = (-1)^{n-1} \prod_{1\le j<n}
    \left(1-\frac{d!}{(dj+1)\cdots(dj+d)}\right)
    \qquad(n\ge1).
\]
It follows that for $n\ge1$
\begin{align}\label{mu-exact}
    \mu_n
    = \sum_{1\le k\le n} \binom{n}{k}\tilde{\mu}_k
    = \sum_{1\le k\le n} \binom{n}{k} (-1)^{k-1}
    \prod_{1\le j<k} \left(1-\frac{d!}
    {(dj+1)\cdots(dj+d)}\right).
\end{align}
This is an identity with exponential cancelation terms; cf.\
\cite{Gnedin07}. In the special case when $d=2$, we have an
identity
\[
    \mu_n = \frac{H_n+2}{3}.
\]
No such simple expression is available for $d\ge3$ since there are
complex-conjugate zeros; see (\ref{d-simplex-mean}).

\paragraph{Exact solution of the general recurrence.}
In general, consider the recurrence
\[
    a_n = b_n + \sum_{0\le k<n} \pi_{n,k} a_k\qquad(n\ge1),
\]
with $a_0=0$. Then the same approach used above leads to the
recurrence
\[
    \tilde{a}_{n+1}
    = -\left(1-\frac{d!}{(dn+1)\cdots(dn+d)}
    \right) \tilde{a}_n + \tilde{b}_n+\tilde{b}_{n+1},
\]
which by iteration gives
\begin{align*}
    \tilde{a}_{n+1}
    &= \sum_{0\le k\le n}(-1)^k
    \left(\tilde{b}_{n-k}+\tilde{b}_{n-k+1}\right)
    \prod_{0\le j<k} \left(1-\frac{d!}
    {(d(n-j)+1)\cdots(d(n-j)+d)}\right),
\end{align*}
by defining $b_0=\tilde{b}_0=0$. Then we obtain the closed-form
solution
\[
    a_n
    = \sum_{1\le k\le n} \binom{n}{k} \tilde{a}_k.
\]
A similar theory of ``d-analogue'' to that presented in
\cite{FLLS95} can be developed (by replacing $2^d/j^d$ there by
$d!/((dj+1)\cdots(dj+d))$).

However, this type of calculations becomes more involved for higher
moments.

\paragraph{Asymptotics of $\mu_n$.} We now look at the asymptotics of
$\mu_n$. To that purpose, we need a better expression for the finite
product in the sum-expression (\ref{mu-exact}).

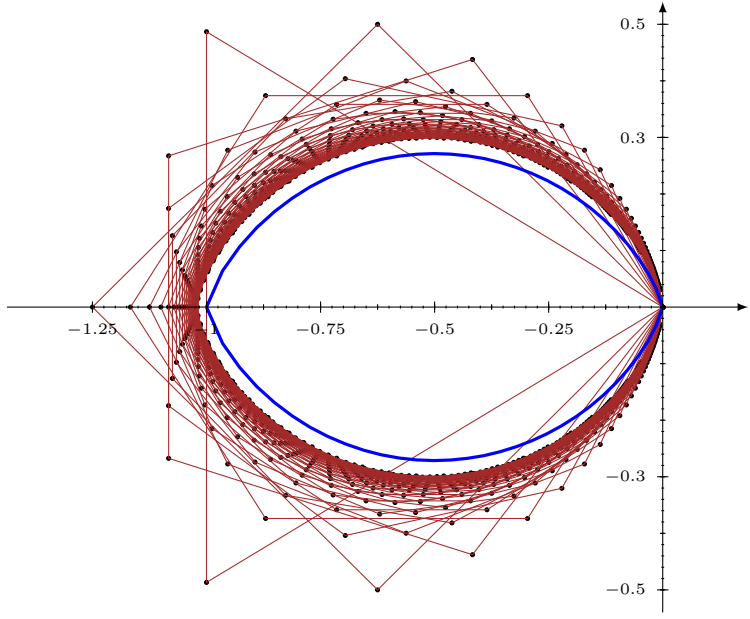
\begin{figure}[!htbp]
\begin{center}
\begin{tikzpicture}[scale=1.5]
\newcommand{\FONTSIZE}{\fontsize{6pt}{\baselineskip}\selectfont}
\definecolor{gen}{rgb}{0,0,0}
\draw[gen] plot[only marks,mark=*,mark size=.5] coordinates{(1.000,
0.066)(1.000, 4.934)(5.000, 2.500)};%
\draw[gen] plot[only marks,mark=*,mark size=.5] coordinates{(0.000,
2.500)(2.500, 0.000)(2.500, 5.000)(5.000, 2.500)};
\draw[gen] plot[only marks,mark=*,mark size=.5] coordinates{(0.668,
1.163)(0.668, 3.837)(3.332, 0.311)(3.332, 4.689)(5.000, 2.500)};
\draw[gen] plot[only marks,mark=*,mark size=.5] coordinates{(0.333,
2.500)(1.518, 0.631)(1.518, 4.369)(3.815, 0.631)(3.815,4.369)
(5.000, 2.500)};%
\draw[gen] plot[only marks,mark=*,mark size=.5] coordinates{(0.667,
1.628)(0.667, 3.372)(2.217, 0.481)(2.217, 4.519)(4.116,0.898)
(4.116, 4.102)(5.000, 2.500)};%
\draw[gen] plot[only marks,mark=*,mark size=.5] coordinates{(0.500,
2.500)(1.186, 1.113)(1.186, 3.887)(2.750, 0.501)(2.750,4.499)
(4.314, 1.113)(4.314, 3.887)(5.000, 2.500)};%
\draw[gen] plot[only marks,mark=*,mark size=.5] coordinates{(0.701,
1.869)(0.701, 3.131)(1.695, 0.836)(1.695, 4.164)(3.152,0.592)
(3.152,4.408)(4.452, 1.286)(4.452, 3.714)(5.000, 2.500)};%
\draw[gen] plot[only marks,mark=*,mark size=.5] coordinates{(0.600,
2.500)(1.049, 1.426)(1.049, 3.574)(2.141, 0.708)(2.141,4.292)
(3.459, 0.708)(3.459, 4.292)(4.551, 1.426)(4.551, 3.574)(5.000,
2.500)};%
\draw[gen] plot[only marks,mark=*,mark size=.5] coordinates{(0.736,
2.012)(0.736, 2.988)(1.425, 1.130)(1.425, 3.870)(2.518,
0.670)(2.518, 4.330)(3.696, 0.829)(3.696, 4.171)(4.624,
1.540)(4.624, 3.460)(5.000, 2.500)};%
\draw[gen] plot[only marks,mark=*,mark size=.5] coordinates{(0.667,
2.500)(0.986, 1.635)(0.986, 3.365)(1.784, 0.945)(1.784,
4.055)(2.833, 0.684)(2.833, 4.316)(3.882, 0.945)(3.882,
4.055)(4.681, 1.635)(4.681, 3.365)(5.000, 2.500)};%
\draw[gen] plot[only marks,mark=*,mark size=.5] coordinates{(0.765,
2.106)(0.765, 2.894)(1.273, 1.353)(1.273, 3.647)(2.110,
0.839)(2.110, 4.161)(3.096, 0.728)(3.096, 4.272)(4.031,
1.053)(4.031, 3.947)(4.725, 1.715)(4.725, 3.285)(5.000, 2.500)};%
\draw[gen] plot[only marks,mark=*,mark size=.5] coordinates{(0.714,
2.500)(0.954, 1.783)(0.954, 3.217)(1.563, 1.152)(1.563,
3.848)(2.399, 0.787)(2.399, 4.213)(3.315, 0.787)(3.315,
4.213)(4.152, 1.152)(4.152, 3.848)(4.760, 1.783)(4.760,
3.217)(5.000, 2.500)};%
\draw[gen] plot[only marks,mark=*,mark size=.5] coordinates{(0.789,
2.172)(0.789, 2.828)(1.181, 1.523)(1.181, 3.477)(1.839,
1.014)(1.839, 3.986)(2.653, 0.772)(2.653, 4.228)(3.500,
0.854)(3.500, 4.146)(4.250, 1.241)(4.250, 3.759)(4.789,
1.841)(4.789, 3.159)(5.000, 2.500)};%
\draw[gen] plot[only marks,mark=*,mark size=.5] coordinates{(0.750,
2.500)(0.938, 1.891)(0.938, 3.109)(1.418, 1.322)(1.418,
3.678)(2.094, 0.924)(2.094, 4.076)(2.875, 0.781)(2.875,
4.219)(3.656, 0.924)(3.656, 4.076)(4.332, 1.322)(4.332,
3.678)(4.812, 1.891)(4.812, 3.109)(5.000, 2.500)};%
\draw[gen] plot[only marks,mark=*,mark size=.5] coordinates{(0.809,
2.221)(0.809, 2.779)(1.121, 1.654)(1.121, 3.346)(1.652,
1.171)(1.652, 3.829)(2.328, 0.870)(2.328, 4.130)(3.069,
0.807)(3.069, 4.193)(3.788, 0.993)(3.788, 4.007)(4.401,
1.395)(4.401, 3.605)(4.832, 1.935)(4.832, 3.065)(5.000, 2.500)};
\draw[gen] plot[only marks,mark=*,mark size=.5] coordinates{(0.778,
2.500)(0.929, 1.973)(0.929, 3.027)(1.318, 1.461)(1.318, 3.539)
(1.875, 1.061)(1.875, 3.939)(2.539, 0.843)(2.539, 4.157)
(3.239,0.843)(3.239, 4.157)(3.902, 1.061)(3.902, 3.939)
(4.459,1.461)(4.459, 3.539)(4.848, 1.973)(4.848, 3.027)
(5.000,2.500)};%
\draw[gen] plot[only marks,mark=*,mark size=.5] coordinates{(0.826,
2.257)(0.826, 2.743)(1.081, 1.758)(1.081, 3.242)(1.518,
1.308)(1.518, 3.692)(2.086, 0.984)(2.086, 4.016)(2.729,
0.836)(2.729, 4.164)(3.388, 0.886)(3.388, 4.114)(4.001,
1.126)(4.001, 3.874)(4.509, 1.521)(4.509, 3.479)(4.862,
2.007)(4.862, 2.993)(5.000, 2.500)}; %
\draw[gen] plot[only marks,mark=*,mark size=.5] coordinates{(0.800,
2.500)(0.925, 2.038)(0.925, 2.962)(1.248, 1.575)(1.248,
3.425)(1.714, 1.188)(1.714, 3.812)(2.281, 0.932)(2.281,
4.068)(2.900, 0.843)(2.900, 4.157)(3.519, 0.932)(3.519,
4.068)(4.086, 1.188)(4.086, 3.812)(4.552, 1.575)(4.552,
3.425)(4.875, 2.038)(4.875, 2.962)(5.000, 2.500)};%
\draw[gen] plot[only marks,mark=*,mark size=.5] coordinates{(0.840,
2.286)(0.840, 2.714)(1.053, 1.841)(1.053, 3.159)(1.420,
1.424)(1.420, 3.576)(1.902, 1.096)(1.902, 3.904)(2.461,
0.900)(2.461, 4.100)(3.054, 0.860)(3.054, 4.140)(3.635,
0.980)(3.635, 4.020)(4.161, 1.246)(4.161, 3.754)(4.590,
1.625)(4.590, 3.375)(4.885, 2.065)(4.885, 2.935)(5.000, 2.500)};%
\draw[gen] plot[only marks,mark=*,mark size=.5] coordinates{(0.818,
2.500)(0.924, 2.089)(0.924, 2.911)(1.196, 1.670)(1.196,
3.330)(1.592, 1.300)(1.592, 3.700)(2.080, 1.028)(2.080,
3.972)(2.626, 0.884)(2.626, 4.116)(3.192, 0.884)(3.192,
4.116)(3.738, 1.028)(3.738, 3.972)(4.226, 1.300)(4.226,
3.700)(4.622, 1.670)(4.622, 3.330)(4.894, 2.089)(4.894,
2.911)(5.000, 2.500)};%
\draw[gen] plot[only marks,mark=*,mark size=.5] coordinates{(0.852,
2.309)(0.852, 2.691)(1.033, 1.909)(1.033, 3.091)(1.346,
1.523)(1.346, 3.477)(1.760, 1.201)(1.760, 3.799)(2.248,
0.979)(2.248, 4.021)(2.778, 0.880)(2.778, 4.120)(3.317,
0.913)(3.317, 4.087)(3.830, 1.076)(3.830, 3.924)(4.284,
1.352)(4.284, 3.648)(4.651, 1.711)(4.651, 3.289)(4.902,
2.111)(4.902, 2.889)(5.000, 2.500)};%
\draw[gen] plot[only marks,mark=*,mark size=.5] coordinates{(0.833,
2.500)(0.924, 2.131)(0.924, 2.869)(1.157, 1.749)(1.157,
3.251)(1.498, 1.400)(1.498, 3.600)(1.922, 1.123)(1.922,
3.877)(2.404, 0.946)(2.404, 4.054)(2.917, 0.885)(2.917,
4.115)(3.429, 0.946)(3.429, 4.054)(3.911, 1.123)(3.911,
3.877)(4.336, 1.400)(4.336, 3.600)(4.676, 1.749)(4.676,
3.251)(4.909, 2.131)(4.909, 2.869)(5.000, 2.500)};%
\draw[gen] plot[only marks,mark=*,mark size=.5] coordinates{(0.862,
2.328)(0.862, 2.672)(1.019, 1.966)(1.019, 3.034)(1.289,
1.608)(1.289, 3.392)(1.648, 1.297)(1.648, 3.703)(2.077,
1.062)(2.077, 3.938)(2.550, 0.925)(2.550, 4.075)(3.044,
0.897)(3.044, 4.103)(3.531, 0.981)(3.531, 4.019)(3.985,
1.169)(3.985, 3.831)(4.382, 1.445)(4.382, 3.555)(4.699,
1.783)(4.699, 3.217)(4.916, 2.149)(4.916, 2.851)(5.000, 2.500)};%
\draw[gen] plot[only marks,mark=*,mark size=.5] coordinates{(0.846,
2.500)(0.925, 2.166)(0.925, 2.834)(1.127, 1.815)(1.127,
3.185)(1.424, 1.487)(1.424, 3.513)(1.795, 1.213)(1.795,
3.787)(2.223, 1.017)(2.223, 3.983)(2.686, 0.915)(2.686,
4.085)(3.160, 0.915)(3.160, 4.085)(3.623, 1.017)(3.623,
3.983)(4.051, 1.213)(4.051, 3.787)(4.423, 1.487)(4.423,
3.513)(4.719, 1.815)(4.719, 3.185)(4.921, 2.166)(4.921,
2.834)(5.000, 2.500)};%
\draw[gen] plot[only marks,mark=*,mark size=.5] coordinates{(0.871,
2.344)(0.871, 2.656)(1.008, 2.013)(1.008, 2.987)(1.244,
1.682)(1.244, 3.318)(1.559, 1.383)(1.559, 3.617)(1.938,
1.144)(1.938, 3.856)(2.362, 0.983)(2.362, 4.017)(2.812,
0.912)(2.812, 4.088)(3.267, 0.936)(3.267, 4.064)(3.707,
1.053)(3.707, 3.947)(4.110, 1.255)(4.110, 3.745)(4.459,
1.526)(4.459, 3.474)(4.737, 1.845)(4.737, 3.155)(4.926,
2.181)(4.926, 2.819)(5.000, 2.500)};%
\draw[gen] plot[only marks,mark=*,mark size=.5] coordinates{(0.857,
2.500)(0.926, 2.195)(0.926, 2.805)(1.104, 1.872)(1.104,
3.128)(1.364, 1.564)(1.364, 3.436)(1.693, 1.296)(1.693,
3.704)(2.074, 1.090)(2.074, 3.910)(2.492, 0.960)(2.492,
4.040)(2.929, 0.916)(2.929, 4.084)(3.365, 0.960)(3.365,
4.040)(3.783, 1.090)(3.783, 3.910)(4.164, 1.296)(4.164,
3.704)(4.493, 1.564)(4.493, 3.436)(4.753, 1.872)(4.753,
3.128)(4.931, 2.195)(4.931, 2.805)(5.000, 2.500)};%
\draw[gen] plot[only marks,mark=*,mark size=.5] coordinates{(0.879,
2.357)(0.879, 2.643)(1.000, 2.054)(1.000, 2.946)(1.208,
1.745)(1.208, 3.255)(1.487, 1.461)(1.487, 3.539)(1.823,
1.223)(1.823, 3.777)(2.204, 1.047)(2.204, 3.953)(2.614,
0.946)(2.614, 4.054)(3.037, 0.925)(3.037, 4.075)(3.455,
0.987)(3.455, 4.013)(3.853, 1.126)(3.853, 3.874)(4.214,
1.335)(4.214, 3.665)(4.523, 1.598)(4.523, 3.402)(4.768,
1.898)(4.768, 3.102)(4.935, 2.208)(4.935, 2.792)(5.000, 2.500)};
\draw[gen] plot[only marks,mark=*,mark size=.5] coordinates{(0.867,
2.500)(0.928, 2.220)(0.928, 2.780)(1.086, 1.922)(1.086,
3.078)(1.316, 1.631)(1.316, 3.369)(1.608, 1.372)(1.608,
3.628)(1.950, 1.162)(1.950, 3.838)(2.328, 1.014)(2.328,
3.986)(2.729, 0.938)(2.729, 4.062)(3.137, 0.938)(3.137,
4.062)(3.538, 1.014)(3.538, 3.986)(3.916, 1.162)(3.916,
3.838)(4.258, 1.372)(4.258, 3.628)(4.550, 1.631)(4.550,
3.369)(4.781, 1.922)(4.781, 3.078)(4.938, 2.220)(4.938,
2.780)(5.000, 2.500)};%
\draw[gen] plot[only marks,mark=*,mark size=.5] coordinates{(0.886,
2.369)(0.886, 2.631)(0.994, 2.088)(0.994, 2.912)(1.179,
1.801)(1.179, 3.199)(1.427, 1.530)(1.427, 3.470)(1.729,
1.296)(1.729, 3.704)(2.072, 1.112)(2.072, 3.888)(2.446,
0.990)(2.446, 4.010)(2.837, 0.937)(2.837, 4.063)(3.231,
0.955)(3.231, 4.045)(3.615, 1.043)(3.615, 3.957)(3.975,
1.197)(3.975, 3.803)(4.299, 1.407)(4.299, 3.593)(4.575,
1.662)(4.575, 3.338)(4.793, 1.944)(4.793, 3.056)(4.942,
2.231)(4.942, 2.769)(5.000, 2.500)};%
\draw[gen] plot[only marks,mark=*,mark size=.5] coordinates{(0.875,
2.500)(0.930, 2.242)(0.930, 2.758)(1.071, 1.964)(1.071,
3.036)(1.277, 1.691)(1.277, 3.309)(1.539, 1.441)(1.539,
3.559)(1.846, 1.231)(1.846, 3.769)(2.190, 1.072)(2.190,
3.928)(2.557, 0.974)(2.557, 4.026)(2.938, 0.940)(2.938,
4.060)(3.318, 0.974)(3.318, 4.026)(3.685, 1.072)(3.685,
3.928)(4.029, 1.231)(4.029, 3.769)(4.336, 1.441)(4.336,
3.559)(4.598, 1.691)(4.598, 3.309)(4.804, 1.964)(4.804,
3.036)(4.945, 2.242)(4.945, 2.758)(5.000, 2.500)};%
\draw[gen] plot[only marks,mark=*,mark size=.5] coordinates{(0.893,
2.378)(0.893, 2.622)(0.989, 2.118)(0.989, 2.882)(1.155,
1.849)(1.155, 3.151)(1.378, 1.592)(1.378, 3.408)(1.649,
1.364)(1.649, 3.636)(1.961, 1.176)(1.961, 3.824)(2.302,
1.041)(2.302, 3.959)(2.663, 0.963)(2.663, 4.037)(3.032,
0.947)(3.032, 4.053)(3.398, 0.994)(3.398, 4.006)(3.751,
1.102)(3.751, 3.898)(4.078, 1.264)(4.078, 3.736)(4.371,
1.473)(4.371, 3.527)(4.619, 1.718)(4.619, 3.282)(4.815,
1.984)(4.815, 3.016)(4.948, 2.251)(4.948, 2.749)(5.000, 2.500)};
\draw[gen] plot[only marks,mark=*,mark size=.5] coordinates{(0.882,
2.500)(0.932, 2.260)(0.932, 2.740)(1.059, 2.002)(1.059,
2.998)(1.244, 1.744)(1.244, 3.256)(1.480, 1.504)(1.480,
3.496)(1.758, 1.296)(1.758, 3.704)(2.071, 1.131)(2.071,
3.869)(2.409, 1.016)(2.409, 3.984)(2.762, 0.958)(2.762,
4.042)(3.120, 0.958)(3.120, 4.042)(3.473, 1.016)(3.473,
3.984)(3.811, 1.131)(3.811, 3.869)(4.124, 1.296)(4.124,
3.704)(4.402, 1.504)(4.402, 3.496)(4.638, 1.744)(4.638,
3.256)(4.824, 2.002)(4.824, 2.998)(4.950, 2.260)(4.950,
2.740)(5.000, 2.500)};%
\draw[gen] plot[only marks,mark=*,mark size=.5] coordinates{(0.898,
2.387)(0.898, 2.613)(0.985, 2.145)(0.985, 2.855)(1.135,
1.892)(1.135, 3.108)(1.336, 1.648)(1.336, 3.352)(1.582,
1.426)(1.582, 3.574)(1.865, 1.238)(1.865, 3.762)(2.178,
1.094)(2.178, 3.906)(2.511, 0.998)(2.511, 4.002)(2.856,
0.957)(2.856, 4.043)(3.203, 0.970)(3.203, 4.030)(3.543,
1.039)(3.543, 3.961)(3.867, 1.160)(3.867, 3.840)(4.166,
1.327)(4.166, 3.673)(4.432, 1.533)(4.432, 3.467)(4.656,
1.768)(4.656, 3.232)(4.832, 2.019)(4.832, 2.981)(4.953,
2.269)(4.953, 2.731)(5.000, 2.500)};%
\draw[gen] plot[only marks,mark=*,mark size=.5] coordinates{(0.889,
2.500)(0.934, 2.276)(0.934, 2.724)(1.049, 2.035)(1.049,
2.965)(1.217, 1.791)(1.217, 3.209)(1.430, 1.561)(1.430,
3.439)(1.684, 1.357)(1.684, 3.643)(1.969, 1.189)(1.969,
3.811)(2.280, 1.063)(2.280, 3.937)(2.608, 0.985)(2.608,
4.015)(2.944, 0.959)(2.944, 4.041)(3.281, 0.985)(3.281,
4.015)(3.609, 1.063)(3.609, 3.937)(3.920, 1.189)(3.920,
3.811)(4.205, 1.357)(4.205, 3.643)(4.459, 1.561)(4.459,
3.439)(4.672, 1.791)(4.672, 3.209)(4.840, 2.035)(4.840,
2.965)(4.955, 2.276)(4.955, 2.724)(5.000, 2.500)};%
\draw[gen] plot[only marks,mark=*,mark size=.5] coordinates{(0.903,
2.394)(0.903, 2.606)(0.983, 2.168)(0.983, 2.832)(1.118,
1.930)(1.118, 3.070)(1.301, 1.698)(1.301, 3.302)(1.525,
1.484)(1.525, 3.516)(1.783, 1.297)(1.783, 3.703)(2.070,
1.147)(2.070, 3.853)(2.378, 1.039)(2.378, 3.961)(2.700,
0.977)(2.700, 4.023)(3.028, 0.965)(3.028, 4.035)(3.354,
1.002)(3.354, 3.998)(3.670, 1.087)(3.670, 3.913)(3.968,
1.217)(3.968, 3.783)(4.242, 1.386)(4.242, 3.614)(4.484,
1.588)(4.484, 3.412)(4.688, 1.813)(4.688, 3.187)(4.848,
2.050)(4.848, 2.950)(4.957, 2.284)(4.957, 2.716)(5.000, 2.500)};
\draw[gen] plot[only marks,mark=*,mark size=.5] coordinates{(0.895,
2.500)(0.936, 2.291)(0.936, 2.709)(1.040, 2.064)(1.040,
2.936)(1.193, 1.833)(1.193, 3.167)(1.388, 1.613)(1.388,
3.387)(1.619, 1.414)(1.619, 3.586)(1.881, 1.245)(1.881,
3.755)(2.168, 1.112)(2.168, 3.888)(2.472, 1.020)(2.472,
3.980)(2.788, 0.973)(2.788, 4.027)(3.107, 0.973)(3.107,
4.027)(3.422, 1.020)(3.422, 3.980)(3.727, 1.112)(3.727,
3.888)(4.014, 1.245)(4.014, 3.755)(4.276, 1.414)(4.276,
3.586)(4.507, 1.613)(4.507, 3.387)(4.702, 1.833)(4.702,
3.167)(4.854, 2.064)(4.854, 2.936)(4.959, 2.291)(4.959,
2.709)(5.000, 2.500)};%
\draw[gen] plot[only marks,mark=*,mark size=.5] coordinates{(0.908,
2.401)(0.908, 2.599)(0.980, 2.189)(0.980, 2.811)(1.104,
1.965)(1.104, 3.035)(1.271, 1.744)(1.271, 3.256)(1.476,
1.537)(1.476, 3.463)(1.713, 1.353)(1.713, 3.647)(1.977,
1.199)(1.977, 3.801)(2.262, 1.082)(2.262, 3.918)(2.562,
1.006)(2.562, 3.994)(2.871, 0.973)(2.871, 4.027)(3.181,
0.984)(3.181, 4.016)(3.487, 1.039)(3.487, 3.961)(3.780,
1.136)(3.780, 3.864)(4.056, 1.272)(4.056, 3.728)(4.307,
1.441)(4.307, 3.559)(4.528, 1.638)(4.528, 3.362)(4.715,
1.853)(4.715, 3.147)(4.861, 2.077)(4.861, 2.923)(4.961,
2.297)(4.961, 2.703)(5.000, 2.500)};%
\draw[gen] plot[only marks,mark=*,mark size=.5] coordinates{(0.900,
2.500)(0.938, 2.304)(0.938, 2.696)(1.033, 2.090)(1.033,
2.910)(1.173, 1.872)(1.173, 3.128)(1.352, 1.661)(1.352,
3.339)(1.564, 1.467)(1.564, 3.533)(1.805, 1.298)(1.805,
3.702)(2.069, 1.160)(2.069, 3.840)(2.353, 1.058)(2.353,
3.942)(2.648, 0.996)(2.648, 4.004)(2.950, 0.975)(2.950,
4.025)(3.252, 0.996)(3.252, 4.004)(3.547, 1.058)(3.547,
3.942)(3.831, 1.160)(3.831, 3.840)(4.095, 1.298)(4.095,
3.702)(4.336, 1.467)(4.336, 3.533)(4.548, 1.661)(4.548,
3.339)(4.727, 1.872)(4.727, 3.128)(4.867, 2.090)(4.867,
2.910)(4.962, 2.304)(4.962, 2.696)(5.000, 2.500)};%
\draw[gen] plot[only marks,mark=*,mark size=.5] coordinates{(0.912,
2.407)(0.912, 2.593)(0.979, 2.207)(0.979, 2.793)(1.092,
1.995)(1.092, 3.005)(1.245, 1.785)(1.245, 3.215)(1.433,
1.585)(1.433, 3.415)(1.651, 1.405)(1.651, 3.595)(1.895,
1.250)(1.895, 3.750)(2.160, 1.127)(2.160, 3.873)(2.440,
1.039)(2.440, 3.961)(2.730, 0.990)(2.730, 4.010)(3.025,
0.980)(3.025, 4.020)(3.318, 1.010)(3.318, 3.990)(3.604,
1.079)(3.604, 3.921)(3.878, 1.185)(3.878, 3.815)(4.132,
1.324)(4.132, 3.676)(4.364, 1.492)(4.364, 3.508)(4.567,
1.683)(4.567, 3.317)(4.738, 1.889)(4.738, 3.111)(4.872,
2.102)(4.872, 2.898)(4.964, 2.309)(4.964, 2.691)(5.000, 2.500)};
\draw[gen] plot[only marks,mark=*,mark size=.5] coordinates{(0.905,
2.500)(0.939, 2.315)(0.939, 2.685)(1.027, 2.113)(1.027,
2.887)(1.156, 1.906)(1.156, 3.094)(1.320, 1.704)(1.320,
3.296)(1.515, 1.516)(1.515, 3.484)(1.738, 1.349)(1.738,
3.651)(1.983, 1.209)(1.983, 3.791)(2.247, 1.099)(2.247,
3.901)(2.524, 1.024)(2.524, 3.976)(2.809, 0.986)(2.809,
4.014)(3.096, 0.986)(3.096, 4.014)(3.381, 1.024)(3.381,
3.976)(3.658, 1.099)(3.658, 3.901)(3.922, 1.209)(3.922,
3.791)(4.167, 1.349)(4.167, 3.651)(4.390, 1.516)(4.390,
3.484)(4.585, 1.704)(4.585, 3.296)(4.749, 1.906)(4.749,
3.094)(4.877, 2.113)(4.877, 2.887)(4.965, 2.315)(4.965,
2.685)(5.000, 2.500)};%
\draw[gen] plot[only marks,mark=*,mark size=.5] coordinates{(0.916,
2.412)(0.916, 2.588)(0.977, 2.224)(0.977, 2.776)(1.082,
2.023)(1.082, 2.977)(1.223, 1.822)(1.223, 3.178)(1.396,
1.630)(1.396, 3.370)(1.597, 1.454)(1.597, 3.546)(1.823,
1.300)(1.823, 3.700)(2.069, 1.172)(2.069, 3.828)(2.331,
1.076)(2.331, 3.924)(2.604, 1.013)(2.604, 3.987)(2.883,
0.986)(2.883, 4.014)(3.164, 0.995)(3.164, 4.005)(3.441,
1.040)(3.441, 3.960)(3.709, 1.120)(3.709, 3.880)(3.963,
1.232)(3.963, 3.768)(4.200, 1.374)(4.200, 3.626)(4.414,
1.539)(4.414, 3.461)(4.601, 1.724)(4.601, 3.276)(4.759,
1.922)(4.759, 3.078)(4.882, 2.124)(4.882, 2.876)(4.967,
2.320)(4.967, 2.680)(5.000, 2.500)};%
\draw[gen] plot[only marks,mark=*,mark size=.5] coordinates{(0.909,
2.500)(0.941, 2.325)(0.941, 2.675)(1.023, 2.134)(1.023,
2.866)(1.141, 1.937)(1.141, 3.063)(1.292, 1.744)(1.292,
3.256)(1.473, 1.562)(1.473, 3.438)(1.679, 1.397)(1.679,
3.603)(1.907, 1.256)(1.907, 3.744)(2.153, 1.141)(2.153,
3.859)(2.412, 1.057)(2.412, 3.943)(2.681, 1.005)(2.681,
3.995)(2.955, 0.988)(2.955, 4.012)(3.228, 1.005)(3.228,
3.995)(3.497, 1.057)(3.497, 3.943)(3.757, 1.141)(3.757,
3.859)(4.002, 1.256)(4.002, 3.744)(4.230, 1.397)(4.230,
3.603)(4.436, 1.562)(4.436, 3.438)(4.617, 1.744)(4.617,
3.256)(4.768, 1.937)(4.768, 3.063)(4.887, 2.134)(4.887,
2.866)(4.968, 2.325)(4.968, 2.675)(5.000, 2.500)};%
\draw[gen] plot[only marks,mark=*,mark size=.5] coordinates{(0.919,
2.417)(0.919, 2.583)(0.976, 2.239)(0.976, 2.761)(1.073,
2.048)(1.073, 2.952)(1.204, 1.856)(1.204, 3.144)(1.363,
1.671)(1.363, 3.329)(1.550, 1.499)(1.550, 3.501)(1.759,
1.346)(1.759, 3.654)(1.989, 1.217)(1.989, 3.783)(2.234,
1.114)(2.234, 3.886)(2.490, 1.041)(2.490, 3.959)(2.755,
1.000)(2.755, 4.000)(3.023, 0.992)(3.023, 4.008)(3.289,
1.017)(3.289, 3.983)(3.550, 1.074)(3.550, 3.926)(3.802,
1.162)(3.802, 3.838)(4.039, 1.278)(4.039, 3.722)(4.259,
1.420)(4.259, 3.580)(4.458, 1.583)(4.458, 3.417)(4.631,
1.762)(4.631, 3.238)(4.777, 1.952)(4.777, 3.048)(4.891,
2.144)(4.891, 2.856)(4.969, 2.330)(4.969, 2.670)(5.000, 2.500)};
\draw[gen] plot[only marks,mark=*,mark size=.5] coordinates{(0.913,
2.500)(0.943, 2.335)(0.943, 2.665)(1.018, 2.153)(1.018,
2.847)(1.128, 1.966)(1.128, 3.034)(1.268, 1.780)(1.268,
3.220)(1.435, 1.604)(1.435, 3.396)(1.627, 1.443)(1.627,
3.557)(1.839, 1.301)(1.839, 3.699)(2.069, 1.183)(2.069,
3.817)(2.312, 1.091)(2.312, 3.909)(2.566, 1.029)(2.566,
3.971)(2.826, 0.997)(2.826, 4.003)(3.087, 0.997)(3.087,
4.003)(3.347, 1.029)(3.347, 3.971)(3.601, 1.091)(3.601,
3.909)(3.844, 1.183)(3.844, 3.817)(4.074, 1.301)(4.074,
3.699)(4.286, 1.443)(4.286, 3.557)(4.478, 1.604)(4.478,
3.396)(4.645, 1.780)(4.645, 3.220)(4.785, 1.966)(4.785,
3.034)(4.895, 2.153)(4.895, 2.847)(4.970, 2.335)(4.970,
2.665)(5.000, 2.500)};%
\draw[gen] plot[only marks,mark=*,mark size=.5] coordinates{(0.923,
2.421)(0.923, 2.579)(0.975, 2.252)(0.975, 2.748)(1.065,
2.071)(1.065, 2.929)(1.186, 1.887)(1.186, 3.113)(1.335,
1.709)(1.335, 3.291)(1.508, 1.542)(1.508, 3.458)(1.703,
1.391)(1.703, 3.609)(1.917, 1.260)(1.917, 3.740)(2.146,
1.153)(2.146, 3.847)(2.388, 1.072)(2.388, 3.928)(2.638,
1.020)(2.638, 3.980)(2.893, 0.997)(2.893, 4.003)(3.149,
1.005)(3.149, 3.995)(3.402, 1.042)(3.402, 3.958)(3.649,
1.109)(3.649, 3.891)(3.885, 1.203)(3.885, 3.797)(4.107,
1.323)(4.107, 3.677)(4.312, 1.464)(4.312, 3.536)(4.496,
1.624)(4.496, 3.376)(4.658, 1.797)(4.658, 3.203)(4.793,
1.979)(4.793, 3.021)(4.899, 2.162)(4.899, 2.838)(4.971,
2.339)(4.971, 2.661)(5.000, 2.500)};%
\draw[gen] plot[only marks,mark=*,mark size=.5] coordinates{(0.917,
2.500)(0.944, 2.343)(0.944, 2.657)(1.015, 2.171)(1.015,
2.829)(1.117, 1.992)(1.117, 3.008)(1.247, 1.814)(1.247,
3.186)(1.402, 1.643)(1.402, 3.357)(1.581, 1.485)(1.581,
3.515)(1.779, 1.344)(1.779, 3.656)(1.993, 1.224)(1.993,
3.776)(2.222, 1.127)(2.222, 3.873)(2.462, 1.056)(2.462,
3.944)(2.708, 1.013)(2.708, 3.987)(2.958, 0.999)(2.958,
4.001)(3.208, 1.013)(3.208, 3.987)(3.455, 1.056)(3.455,
3.944)(3.694, 1.127)(3.694, 3.873)(3.923, 1.224)(3.923,
3.776)(4.138, 1.344)(4.138, 3.656)(4.336, 1.485)(4.336,
3.515)(4.514, 1.643)(4.514, 3.357)(4.670, 1.814)(4.670,
3.186)(4.800, 1.992)(4.800, 3.008)(4.902, 2.171)(4.902,
2.829)(4.972, 2.343)(4.972, 2.657)(5.000, 2.500)};
\draw[gen] plot[only marks,mark=*,mark size=.5] coordinates{(0.925,
2.425)(0.925, 2.575)(0.975, 2.264)(0.975, 2.736)(1.059,
2.092)(1.059, 2.908)(1.171, 1.916)(1.171, 3.084)(1.309,
1.744)(1.309, 3.256)(1.471, 1.582)(1.471, 3.418)(1.653,
1.433)(1.653, 3.567)(1.853, 1.302)(1.853, 3.698)(2.068,
1.192)(2.068, 3.808)(2.296, 1.105)(2.296, 3.895)(2.533,
1.044)(2.533, 3.956)(2.775, 1.009)(2.775, 3.991)(3.021,
1.002)(3.021, 3.998)(3.265, 1.023)(3.265, 3.977)(3.505,
1.071)(3.505, 3.929)(3.738, 1.145)(3.738, 3.855)(3.960,
1.244)(3.960, 3.756)(4.168, 1.365)(4.168, 3.635)(4.359,
1.506)(4.359, 3.494)(4.531, 1.662)(4.531, 3.338)(4.681,
1.829)(4.681, 3.171)(4.807, 2.004)(4.807, 2.996)(4.905,
2.179)(4.905, 2.821)(4.973, 2.347)(4.973, 2.653)(5.000, 2.500)};
\draw[gen] plot[only marks,mark=*,mark size=.5] coordinates{(0.920,
2.500)(0.946, 2.351)(0.946, 2.649)(1.011, 2.187)(1.011,
2.813)(1.107, 2.015)(1.107, 2.985)(1.228, 1.845)(1.228,
3.155)(1.373, 1.680)(1.373, 3.320)(1.539, 1.525)(1.539,
3.475)(1.724, 1.386)(1.724, 3.614)(1.926, 1.264)(1.926,
3.736)(2.141, 1.164)(2.141, 3.836)(2.367, 1.086)(2.367,
3.914)(2.601, 1.034)(2.601, 3.966)(2.840, 1.007)(2.840,
3.993)(3.080, 1.007)(3.080, 3.993)(3.319, 1.034)(3.319,
3.966)(3.553, 1.086)(3.553, 3.914)(3.779, 1.164)(3.779,
3.836)(3.994, 1.264)(3.994, 3.736)(4.196, 1.386)(4.196,
3.614)(4.381, 1.525)(4.381, 3.475)(4.547, 1.680)(4.547,
3.320)(4.692, 1.845)(4.692, 3.155)(4.813, 2.015)(4.813,
2.985)(4.909, 2.187)(4.909, 2.813)(4.974, 2.351)(4.974,
2.649)(5.000, 2.500)};%
\definecolor{gen}{rgb}{0.647059,0.164706,0.164706}
\draw[gen] plot coordinates{(1.000, 0.066)(1.000, 4.934)(5.000,
2.500)(1.000, 0.066)};%
\draw[gen] plot coordinates{(0.668, 1.163)(0.668, 3.837)(3.332,
4.689)(5.000, 2.500)(3.332, 0.311)(0.668, 1.163)};%
\draw[gen] plot coordinates{(0.667, 1.628)(0.667, 3.372)(2.217,
4.519)(4.116, 4.102)(5.000, 2.500)(4.116, 0.898)(2.217,
0.481)(0.667, 1.628)};%
\draw[gen] plot coordinates{(0.701, 1.869)(0.701, 3.131)(1.695,
4.164)(3.152, 4.408)(4.452, 3.714)(5.000, 2.500)(4.452,
1.286)(3.152, 0.592)(1.695, 0.836)(0.701, 1.869)};%
\draw[gen] plot coordinates{(0.736, 2.012)(0.736, 2.988)(1.425,
3.870)(2.518, 4.330)(3.696, 4.171)(4.624, 3.460)(5.000,
2.500)(4.624, 1.540)(3.696, 0.829)(2.518, 0.670)(1.425,
1.130)(0.736, 2.012)};%
\draw[gen] plot coordinates{(0.765, 2.106)(0.765, 2.894)(1.273,
3.647)(2.110, 4.161)(3.096, 4.272)(4.031, 3.947)(4.725,
3.285)(5.000, 2.500)(4.725, 1.715)(4.031, 1.053)(3.096,
0.728)(2.110, 0.839)(1.273, 1.353)(0.765, 2.106)};%
\draw[gen] plot coordinates{(0.789, 2.172)(0.789, 2.828)(1.181,
3.477)(1.839, 3.986)(2.653, 4.228)(3.500, 4.146)(4.250,
3.759)(4.789, 3.159)(5.000, 2.500)(4.789, 1.841)(4.250,
1.241)(3.500, 0.854)(2.653, 0.772)(1.839, 1.014)(1.181,
1.523)(0.789, 2.172)};%
\draw[gen] plot coordinates{(0.809, 2.221)(0.809, 2.779)(1.121,
3.346)(1.652, 3.829)(2.328, 4.130)(3.069, 4.193)(3.788,
4.007)(4.401, 3.605)(4.832, 3.065)(5.000, 2.500)(4.832,
1.935)(4.401, 1.395)(3.788, 0.993)(3.069, 0.807)(2.328,
0.870)(1.652, 1.171)(1.121, 1.654)(0.809, 2.221)};%
\draw[gen] plot coordinates{(0.826, 2.257)(0.826, 2.743)(1.081,
3.242)(1.518, 3.692)(2.086, 4.016)(2.729, 4.164)(3.388,
4.114)(4.001, 3.874)(4.509, 3.479)(4.862, 2.993)(5.000,
2.500)(4.862, 2.007)(4.509, 1.521)(4.001, 1.126)(3.388,
0.886)(2.729, 0.836)(2.086, 0.984)(1.518, 1.308)(1.081,
1.758)(0.826, 2.257)};%
\draw[gen] plot coordinates{(0.840, 2.286)(0.840, 2.714)(1.053,
3.159)(1.420, 3.576)(1.902, 3.904)(2.461, 4.100)(3.054,
4.140)(3.635, 4.020)(4.161, 3.754)(4.590, 3.375)(4.885,
2.935)(5.000, 2.500)(4.885, 2.065)(4.590, 1.625)(4.161,
1.246)(3.635, 0.980)(3.054, 0.860)(2.461, 0.900)(1.902,
1.096)(1.420, 1.424)(1.053, 1.841)(0.840, 2.286)};%
\draw[gen] plot coordinates{(0.852, 2.309)(0.852, 2.691)(1.033,
3.091)(1.346, 3.477)(1.760, 3.799)(2.248, 4.021)(2.778,
4.120)(3.317, 4.087)(3.830, 3.924)(4.284, 3.648)(4.651,
3.289)(4.902, 2.889)(5.000, 2.500)(4.902, 2.111)(4.651,
1.711)(4.284, 1.352)(3.830, 1.076)(3.317, 0.913)(2.778,
0.880)(2.248, 0.979)(1.760, 1.201)(1.346, 1.523)(1.033,
1.909)(0.852, 2.309)};%
\draw[gen] plot coordinates{(0.862, 2.328)(0.862, 2.672)(1.019,
3.034)(1.289, 3.392)(1.648, 3.703)(2.077, 3.938)(2.550,
4.075)(3.044, 4.103)(3.531, 4.019)(3.985, 3.831)(4.382,
3.555)(4.699, 3.217)(4.916, 2.851)(5.000, 2.500)(4.916,
2.149)(4.699, 1.783)(4.382, 1.445)(3.985, 1.169)(3.531,
0.981)(3.044, 0.897)(2.550, 0.925)(2.077, 1.062)(1.648,
1.297)(1.289, 1.608)(1.019, 1.966)(0.862, 2.328)};%
\draw[gen] plot coordinates{(0.871, 2.344)(0.871, 2.656)(1.008,
2.987)(1.244, 3.318)(1.559, 3.617)(1.938, 3.856)(2.362,
4.017)(2.812, 4.088)(3.267, 4.064)(3.707, 3.947)(4.110,
3.745)(4.459, 3.474)(4.737, 3.155)(4.926, 2.819)(5.000,
2.500)(4.926, 2.181)(4.737, 1.845)(4.459, 1.526)(4.110,
1.255)(3.707, 1.053)(3.267, 0.936)(2.812, 0.912)(2.362,
0.983)(1.938, 1.144)(1.559, 1.383)(1.244, 1.682)(1.008,
2.013)(0.871, 2.344)};%
\draw[gen] plot coordinates{(0.879, 2.357)(0.879, 2.643)(1.000,
2.946)(1.208, 3.255)(1.487, 3.539)(1.823, 3.777)(2.204,
3.953)(2.614, 4.054)(3.037, 4.075)(3.455, 4.013)(3.853,
3.874)(4.214, 3.665)(4.523, 3.402)(4.768, 3.102)(4.935,
2.792)(5.000, 2.500)(4.935, 2.208)(4.768, 1.898)(4.523,
1.598)(4.214, 1.335)(3.853, 1.126)(3.455, 0.987)(3.037,
0.925)(2.614, 0.946)(2.204, 1.047)(1.823, 1.223)(1.487,
1.461)(1.208, 1.745)(1.000, 2.054)(0.879, 2.357)}; %
\draw[gen] plot coordinates{(0.886, 2.369)(0.886, 2.631)(0.994,
2.912)(1.179, 3.199)(1.427, 3.470)(1.729, 3.704)(2.072,
3.888)(2.446, 4.010)(2.837, 4.063)(3.231, 4.045)(3.615,
3.957)(3.975, 3.803)(4.299, 3.593)(4.575, 3.338)(4.793,
3.056)(4.942, 2.769)(5.000, 2.500)(4.942, 2.231)(4.793,
1.944)(4.575, 1.662)(4.299, 1.407)(3.975, 1.197)(3.615,
1.043)(3.231, 0.955)(2.837, 0.937)(2.446, 0.990)(2.072,
1.112)(1.729, 1.296)(1.427, 1.530)(1.179, 1.801)(0.994,
2.088)(0.886, 2.369)};%
\draw[gen] plot coordinates{(0.893, 2.378)(0.893, 2.622)(0.989,
2.882)(1.155, 3.151)(1.378, 3.408)(1.649, 3.636)(1.961,
3.824)(2.302, 3.959)(2.663, 4.037)(3.032, 4.053)(3.398,
4.006)(3.751, 3.898)(4.078, 3.736)(4.371, 3.527)(4.619,
3.282)(4.815, 3.016)(4.948, 2.749)(5.000, 2.500)(4.948,
2.251)(4.815, 1.984)(4.619, 1.718)(4.371, 1.473)(4.078,
1.264)(3.751, 1.102)(3.398, 0.994)(3.032, 0.947)(2.663,
0.963)(2.302, 1.041)(1.961, 1.176)(1.649, 1.364)(1.378,
1.592)(1.155, 1.849)(0.989, 2.118)(0.893, 2.378)};%
\draw[gen] plot coordinates{(0.898, 2.387)(0.898, 2.613)(0.985,
2.855)(1.135, 3.108)(1.336, 3.352)(1.582, 3.574)(1.865,
3.762)(2.178, 3.906)(2.511, 4.002)(2.856, 4.043)(3.203,
4.030)(3.543, 3.961)(3.867, 3.840)(4.166, 3.673)(4.432,
3.467)(4.656, 3.232)(4.832, 2.981)(4.953, 2.731)(5.000,
2.500)(4.953, 2.269)(4.832, 2.019)(4.656, 1.768)(4.432,
1.533)(4.166, 1.327)(3.867, 1.160)(3.543, 1.039)(3.203,
0.970)(2.856, 0.957)(2.511, 0.998)(2.178, 1.094)(1.865,
1.238)(1.582, 1.426)(1.336, 1.648)(1.135, 1.892)(0.985,
2.145)(0.898, 2.387)};%
\draw[gen] plot coordinates{(0.903, 2.394)(0.903, 2.606)(0.983,
2.832)(1.118, 3.070)(1.301, 3.302)(1.525, 3.516)(1.783,
3.703)(2.070, 3.853)(2.378, 3.961)(2.700, 4.023)(3.028,
4.035)(3.354, 3.998)(3.670, 3.913)(3.968, 3.783)(4.242,
3.614)(4.484, 3.412)(4.688, 3.187)(4.848, 2.950)(4.957,
2.716)(5.000, 2.500)(4.957, 2.284)(4.848, 2.050)(4.688,
1.813)(4.484, 1.588)(4.242, 1.386)(3.968, 1.217)(3.670,
1.087)(3.354, 1.002)(3.028, 0.965)(2.700, 0.977)(2.378,
1.039)(2.070, 1.147)(1.783, 1.297)(1.525, 1.484)(1.301,
1.698)(1.118, 1.930)(0.983, 2.168)(0.903, 2.394)};%
\draw[gen] plot coordinates{(0.908, 2.401)(0.908, 2.599)(0.980,
2.811)(1.104, 3.035)(1.271, 3.256)(1.476, 3.463)(1.713,
3.647)(1.977, 3.801)(2.262, 3.918)(2.562, 3.994)(2.871,
4.027)(3.181, 4.016)(3.487, 3.961)(3.780, 3.864)(4.056,
3.728)(4.307, 3.559)(4.528, 3.362)(4.715, 3.147)(4.861,
2.923)(4.961, 2.703)(5.000, 2.500)(4.961, 2.297)(4.861,
2.077)(4.715, 1.853)(4.528, 1.638)(4.307, 1.441)(4.056,
1.272)(3.780, 1.136)(3.487, 1.039)(3.181, 0.984)(2.871,
0.973)(2.562, 1.006)(2.262, 1.082)(1.977, 1.199)(1.713,
1.353)(1.476, 1.537)(1.271, 1.744)(1.104, 1.965)(0.980,
2.189)(0.908, 2.401)};%
\draw[gen] plot coordinates{(0.912, 2.407)(0.912, 2.593)(0.979,
2.793)(1.092, 3.005)(1.245, 3.215)(1.433, 3.415)(1.651,
3.595)(1.895, 3.750)(2.160, 3.873)(2.440, 3.961)(2.730,
4.010)(3.025, 4.020)(3.318, 3.990)(3.604, 3.921)(3.878,
3.815)(4.132, 3.676)(4.364, 3.508)(4.567, 3.317)(4.738,
3.111)(4.872, 2.898)(4.964, 2.691)(5.000, 2.500)(4.964,
2.309)(4.872, 2.102)(4.738, 1.889)(4.567, 1.683)(4.364,
1.492)(4.132, 1.324)(3.878, 1.185)(3.604, 1.079)(3.318,
1.010)(3.025, 0.980)(2.730, 0.990)(2.440, 1.039)(2.160,
1.127)(1.895, 1.250)(1.651, 1.405)(1.433, 1.585)(1.245,
1.785)(1.092, 1.995)(0.979, 2.207)(0.912, 2.407)};%

\draw[gen] plot coordinates{(0.916, 2.412)(0.916, 2.588)(0.977,
2.776)(1.082, 2.977)(1.223, 3.178)(1.396, 3.370)(1.597,
3.546)(1.823, 3.700)(2.069, 3.828)(2.331, 3.924)(2.604,
3.987)(2.883, 4.014)(3.164, 4.005)(3.441, 3.960)(3.709,
3.880)(3.963, 3.768)(4.200, 3.626)(4.414, 3.461)(4.601,
3.276)(4.759, 3.078)(4.882, 2.876)(4.967, 2.680)(5.000,
2.500)(4.967, 2.320)(4.882, 2.124)(4.759, 1.922)(4.601,
1.724)(4.414, 1.539)(4.200, 1.374)(3.963, 1.232)(3.709,
1.120)(3.441, 1.040)(3.164, 0.995)(2.883, 0.986)(2.604,
1.013)(2.331, 1.076)(2.069, 1.172)(1.823, 1.300)(1.597,
1.454)(1.396, 1.630)(1.223, 1.822)(1.082, 2.023)(0.977,
2.224)(0.916, 2.412)};%
\draw[gen] plot coordinates{(0.919, 2.417)(0.919, 2.583)(0.976,
2.761)(1.073, 2.952)(1.204, 3.144)(1.363, 3.329)(1.550,
3.501)(1.759, 3.654)(1.989, 3.783)(2.234, 3.886)(2.490,
3.959)(2.755, 4.000)(3.023, 4.008)(3.289, 3.983)(3.550,
3.926)(3.802, 3.838)(4.039, 3.722)(4.259, 3.580)(4.458,
3.417)(4.631, 3.238)(4.777, 3.048)(4.891, 2.856)(4.969,
2.670)(5.000, 2.500)(4.969, 2.330)(4.891, 2.144)(4.777,
1.952)(4.631, 1.762)(4.458, 1.583)(4.259, 1.420)(4.039,
1.278)(3.802, 1.162)(3.550, 1.074)(3.289, 1.017)(3.023,
0.992)(2.755, 1.000)(2.490, 1.041)(2.234, 1.114)(1.989,
1.217)(1.759, 1.346)(1.550, 1.499)(1.363, 1.671)(1.204,
1.856)(1.073, 2.048)(0.976, 2.239)(0.919, 2.417)};%
\draw[gen] plot coordinates{(0.923, 2.421)(0.923, 2.579)(0.975,
2.748)(1.065, 2.929)(1.186, 3.113)(1.335, 3.291)(1.508,
3.458)(1.703, 3.609)(1.917, 3.740)(2.146, 3.847)(2.388,
3.928)(2.638, 3.980)(2.893, 4.003)(3.149, 3.995)(3.402,
3.958)(3.649, 3.891)(3.885, 3.797)(4.107, 3.677)(4.312,
3.536)(4.496, 3.376)(4.658, 3.203)(4.793, 3.021)(4.899,
2.838)(4.971, 2.661)(5.000, 2.500)(4.971, 2.339)(4.899,
2.162)(4.793, 1.979)(4.658, 1.797)(4.496, 1.624)(4.312,
1.464)(4.107, 1.323)(3.885, 1.203)(3.649, 1.109)(3.402,
1.042)(3.149, 1.005)(2.893, 0.997)(2.638, 1.020)(2.388,
1.072)(2.146, 1.153)(1.917, 1.260)(1.703, 1.391)(1.508,
1.542)(1.335, 1.709)(1.186, 1.887)(1.065, 2.071)(0.975,
2.252)(0.923, 2.421)};%
\draw[gen] plot coordinates{(0.925, 2.425)(0.925, 2.575)(0.975,
2.736)(1.059, 2.908)(1.171, 3.084)(1.309, 3.256)(1.471,
3.418)(1.653, 3.567)(1.853, 3.698)(2.068, 3.808)(2.296,
3.895)(2.533, 3.956)(2.775, 3.991)(3.021, 3.998)(3.265,
3.977)(3.505, 3.929)(3.738, 3.855)(3.960, 3.756)(4.168,
3.635)(4.359, 3.494)(4.531, 3.338)(4.681, 3.171)(4.807,
2.996)(4.905, 2.821)(4.973, 2.653)(5.000, 2.500)(4.973,
2.347)(4.905, 2.179)(4.807, 2.004)(4.681, 1.829)(4.531,
1.662)(4.359, 1.506)(4.168, 1.365)(3.960, 1.244)(3.738,
1.145)(3.505, 1.071)(3.265, 1.023)(3.021, 1.002)(2.775,
1.009)(2.533, 1.044)(2.296, 1.105)(2.068, 1.192)(1.853,
1.302)(1.653, 1.433)(1.471, 1.582)(1.309, 1.744)(1.171,
1.916)(1.059, 2.092)(0.975, 2.264)(0.925, 2.425)};%
\draw[gen] plot coordinates{(0.000, 2.500)(2.500, 5.000)(5.000,
2.500)(2.500, 0.000)(0.000, 2.500)};%
\draw[gen] plot coordinates{(0.333, 2.500)(1.518, 4.369)(3.815,
4.369)(5.000, 2.500)(3.815, 0.631)(1.518, 0.631)(0.333, 2.500)};
\draw[gen] plot coordinates{(0.500, 2.500)(1.186, 3.887)(2.750,
4.499)(4.314, 3.887)(5.000, 2.500)(4.314, 1.113)(2.750,
0.501)(1.186, 1.113)(0.500, 2.500)};%
\draw[gen] plot coordinates{(0.600, 2.500)(1.049, 3.574)(2.141,
4.292)(3.459, 4.292)(4.551, 3.574)(5.000, 2.500)(4.551,
1.426)(3.459, 0.708)(2.141, 0.708)(1.049, 1.426)(0.600, 2.500)};
\draw[gen] plot coordinates{(0.667, 2.500)(0.986, 3.365)(1.784,
4.055)(2.833, 4.316)(3.882, 4.055)(4.681, 3.365)(5.000,
2.500)(4.681, 1.635)(3.882, 0.945)(2.833, 0.684)(1.784,
0.945)(0.986, 1.635)(0.667, 2.500)};%
\draw[gen] plot coordinates{(0.714, 2.500)(0.954, 3.217)(1.563,
3.848)(2.399, 4.213)(3.315, 4.213)(4.152, 3.848)(4.760,
3.217)(5.000, 2.500)(4.760, 1.783)(4.152, 1.152)(3.315,
0.787)(2.399, 0.787)(1.563, 1.152)(0.954, 1.783)(0.714, 2.500)};
\draw[gen] plot coordinates{(0.750, 2.500)(0.938, 3.109)(1.418,
3.678)(2.094, 4.076)(2.875, 4.219)(3.656, 4.076)(4.332,
3.678)(4.812, 3.109)(5.000, 2.500)(4.812, 1.891)(4.332,
1.322)(3.656, 0.924)(2.875, 0.781)(2.094, 0.924)(1.418,
1.322)(0.938, 1.891)(0.750, 2.500)};%
\draw[gen] plot coordinates{(0.778, 2.500)(0.929, 3.027)(1.318,
3.539)(1.875, 3.939)(2.539, 4.157)(3.239, 4.157)(3.902,
3.939)(4.459, 3.539)(4.848, 3.027)(5.000, 2.500)(4.848,
1.973)(4.459, 1.461)(3.902, 1.061)(3.239, 0.843)(2.539,
0.843)(1.875, 1.061)(1.318, 1.461)(0.929, 1.973)(0.778, 2.500)};
\draw[gen] plot coordinates{(0.800, 2.500)(0.925, 2.962)(1.248,
3.425)(1.714, 3.812)(2.281, 4.068)(2.900, 4.157)(3.519,
4.068)(4.086, 3.812)(4.552, 3.425)(4.875, 2.962)(5.000,
2.500)(4.875, 2.038)(4.552, 1.575)(4.086, 1.188)(3.519,
0.932)(2.900, 0.843)(2.281, 0.932)(1.714, 1.188)(1.248,
1.575)(0.925, 2.038)(0.800, 2.500)};%
\draw[gen] plot coordinates{(0.818, 2.500)(0.924, 2.911)(1.196,
3.330)(1.592, 3.700)(2.080, 3.972)(2.626, 4.116)(3.192,
4.116)(3.738, 3.972)(4.226, 3.700)(4.622, 3.330)(4.894,
2.911)(5.000, 2.500)(4.894, 2.089)(4.622, 1.670)(4.226,
1.300)(3.738, 1.028)(3.192, 0.884)(2.626, 0.884)(2.080,
1.028)(1.592, 1.300)(1.196, 1.670)(0.924, 2.089)(0.818, 2.500)};
\draw[gen] plot coordinates{(0.833, 2.500)(0.924, 2.869)(1.157,
3.251)(1.498, 3.600)(1.922, 3.877)(2.404, 4.054)(2.917,
4.115)(3.429, 4.054)(3.911, 3.877)(4.336, 3.600)(4.676,
3.251)(4.909, 2.869)(5.000, 2.500)(4.909, 2.131)(4.676,
1.749)(4.336, 1.400)(3.911, 1.123)(3.429, 0.946)(2.917,
0.885)(2.404, 0.946)(1.922, 1.123)(1.498, 1.400)(1.157,
1.749)(0.924, 2.131)(0.833, 2.500)};%
\draw[gen] plot coordinates{(0.846, 2.500)(0.925, 2.834)(1.127,
3.185)(1.424, 3.513)(1.795, 3.787)(2.223, 3.983)(2.686,
4.085)(3.160, 4.085)(3.623, 3.983)(4.051, 3.787)(4.423,
3.513)(4.719, 3.185)(4.921, 2.834)(5.000, 2.500)(4.921,
2.166)(4.719, 1.815)(4.423, 1.487)(4.051, 1.213)(3.623,
1.017)(3.160, 0.915)(2.686, 0.915)(2.223, 1.017)(1.795,
1.213)(1.424, 1.487)(1.127, 1.815)(0.925, 2.166)(0.846, 2.500)};
\draw[gen] plot coordinates{(0.857, 2.500)(0.926, 2.805)(1.104,
3.128)(1.364, 3.436)(1.693, 3.704)(2.074, 3.910)(2.492,
4.040)(2.929, 4.084)(3.365, 4.040)(3.783, 3.910)(4.164,
3.704)(4.493, 3.436)(4.753, 3.128)(4.931, 2.805)(5.000,
2.500)(4.931, 2.195)(4.753, 1.872)(4.493, 1.564)(4.164,
1.296)(3.783, 1.090)(3.365, 0.960)(2.929, 0.916)(2.492,
0.960)(2.074, 1.090)(1.693, 1.296)(1.364, 1.564)(1.104,
1.872)(0.926, 2.195)(0.857, 2.500)};%
\draw[gen] plot coordinates{(0.867, 2.500)(0.928, 2.780)(1.086,
3.078)(1.316, 3.369)(1.608, 3.628)(1.950, 3.838)(2.328,
3.986)(2.729, 4.062)(3.137, 4.062)(3.538, 3.986)(3.916,
3.838)(4.258, 3.628)(4.550, 3.369)(4.781, 3.078)(4.938,
2.780)(5.000, 2.500)(4.938, 2.220)(4.781, 1.922)(4.550,
1.631)(4.258, 1.372)(3.916, 1.162)(3.538, 1.014)(3.137,
0.938)(2.729, 0.938)(2.328, 1.014)(1.950, 1.162)(1.608,
1.372)(1.316, 1.631)(1.086, 1.922)(0.928, 2.220)(0.867, 2.500)};
\draw[gen] plot coordinates{(0.875, 2.500)(0.930, 2.758)(1.071,
3.036)(1.277, 3.309)(1.539, 3.559)(1.846, 3.769)(2.190,
3.928)(2.557, 4.026)(2.938, 4.060)(3.318, 4.026)(3.685,
3.928)(4.029, 3.769)(4.336, 3.559)(4.598, 3.309)(4.804,
3.036)(4.945, 2.758)(5.000, 2.500)(4.945, 2.242)(4.804,
1.964)(4.598, 1.691)(4.336, 1.441)(4.029, 1.231)(3.685,
1.072)(3.318, 0.974)(2.938, 0.940)(2.557, 0.974)(2.190,
1.072)(1.846, 1.231)(1.539, 1.441)(1.277, 1.691)(1.071,
1.964)(0.930, 2.242)(0.875, 2.500)};%
\draw[gen] plot coordinates{(0.882, 2.500)(0.932, 2.740)(1.059,
2.998)(1.244, 3.256)(1.480, 3.496)(1.758, 3.704)(2.071,
3.869)(2.409, 3.984)(2.762, 4.042)(3.120, 4.042)(3.473,
3.984)(3.811, 3.869)(4.124, 3.704)(4.402, 3.496)(4.638,
3.256)(4.824, 2.998)(4.950, 2.740)(5.000, 2.500)(4.950,
2.260)(4.824, 2.002)(4.638, 1.744)(4.402, 1.504)(4.124,
1.296)(3.811, 1.131)(3.473, 1.016)(3.120, 0.958)(2.762,
0.958)(2.409, 1.016)(2.071, 1.131)(1.758, 1.296)(1.480,
1.504)(1.244, 1.744)(1.059, 2.002)(0.932, 2.260)(0.882, 2.500)};%
\draw[gen] plot coordinates{(0.889, 2.500)(0.934, 2.724)(1.049,
2.965)(1.217, 3.209)(1.430, 3.439)(1.684, 3.643)(1.969,
3.811)(2.280, 3.937)(2.608, 4.015)(2.944, 4.041)(3.281,
4.015)(3.609, 3.937)(3.920, 3.811)(4.205, 3.643)(4.459,
3.439)(4.672, 3.209)(4.840, 2.965)(4.955, 2.724)(5.000,
2.500)(4.955, 2.276)(4.840, 2.035)(4.672, 1.791)(4.459,
1.561)(4.205, 1.357)(3.920, 1.189)(3.609, 1.063)(3.281,
0.985)(2.944, 0.959)(2.608, 0.985)(2.280, 1.063)(1.969,
1.189)(1.684, 1.357)(1.430, 1.561)(1.217, 1.791)(1.049,
2.035)(0.934, 2.276)(0.889, 2.500)};%
\draw[gen] plot coordinates{(0.895, 2.500)(0.936, 2.709)(1.040,
2.936)(1.193, 3.167)(1.388, 3.387)(1.619, 3.586)(1.881,
3.755)(2.168, 3.888)(2.472, 3.980)(2.788, 4.027)(3.107,
4.027)(3.422, 3.980)(3.727, 3.888)(4.014, 3.755)(4.276,
3.586)(4.507, 3.387)(4.702, 3.167)(4.854, 2.936)(4.959,
2.709)(5.000, 2.500)(4.959, 2.291)(4.854, 2.064)(4.702,
1.833)(4.507, 1.613)(4.276, 1.414)(4.014, 1.245)(3.727,
1.112)(3.422, 1.020)(3.107, 0.973)(2.788, 0.973)(2.472,
1.020)(2.168, 1.112)(1.881, 1.245)(1.619, 1.414)(1.388,
1.613)(1.193, 1.833)(1.040, 2.064)(0.936, 2.291)(0.895, 2.500)};%
\draw[gen] plot coordinates{(0.900, 2.500)(0.938, 2.696)(1.033,
2.910)(1.173, 3.128)(1.352, 3.339)(1.564, 3.533)(1.805,
3.702)(2.069, 3.840)(2.353, 3.942)(2.648, 4.004)(2.950,
4.025)(3.252, 4.004)(3.547, 3.942)(3.831, 3.840)(4.095,
3.702)(4.336, 3.533)(4.548, 3.339)(4.727, 3.128)(4.867,
2.910)(4.962, 2.696)(5.000, 2.500)(4.962, 2.304)(4.867,
2.090)(4.727, 1.872)(4.548, 1.661)(4.336, 1.467)(4.095,
1.298)(3.831, 1.160)(3.547, 1.058)(3.252, 0.996)(2.950,
0.975)(2.648, 0.996)(2.353, 1.058)(2.069, 1.160)(1.805,
1.298)(1.564, 1.467)(1.352, 1.661)(1.173, 1.872)(1.033,
2.090)(0.938, 2.304)(0.900, 2.500)};%
\draw[gen] plot coordinates{(0.905, 2.500)(0.939, 2.685)(1.027,
2.887)(1.156, 3.094)(1.320, 3.296)(1.515, 3.484)(1.738,
3.651)(1.983, 3.791)(2.247, 3.901)(2.524, 3.976)(2.809,
4.014)(3.096, 4.014)(3.381, 3.976)(3.658, 3.901)(3.922,
3.791)(4.167, 3.651)(4.390, 3.484)(4.585, 3.296)(4.749,
3.094)(4.877, 2.887)(4.965, 2.685)(5.000, 2.500)(4.965,
2.315)(4.877, 2.113)(4.749, 1.906)(4.585, 1.704)(4.390,
1.516)(4.167, 1.349)(3.922, 1.209)(3.658, 1.099)(3.381,
1.024)(3.096, 0.986)(2.809, 0.986)(2.524, 1.024)(2.247,
1.099)(1.983, 1.209)(1.738, 1.349)(1.515, 1.516)(1.320,
1.704)(1.156, 1.906)(1.027, 2.113)(0.939, 2.315)(0.905, 2.500)};%
\draw[gen] plot coordinates{(0.909, 2.500)(0.941, 2.675)(1.023,
2.866)(1.141, 3.063)(1.292, 3.256)(1.473, 3.438)(1.679,
3.603)(1.907, 3.744)(2.153, 3.859)(2.412, 3.943)(2.681,
3.995)(2.955, 4.012)(3.228, 3.995)(3.497, 3.943)(3.757,
3.859)(4.002, 3.744)(4.230, 3.603)(4.436, 3.438)(4.617,
3.256)(4.768, 3.063)(4.887, 2.866)(4.968, 2.675)(5.000,
2.500)(4.968, 2.325)(4.887, 2.134)(4.768, 1.937)(4.617,
1.744)(4.436, 1.562)(4.230, 1.397)(4.002, 1.256)(3.757,
1.141)(3.497, 1.057)(3.228, 1.005)(2.955, 0.988)(2.681,
1.005)(2.412, 1.057)(2.153, 1.141)(1.907, 1.256)(1.679,
1.397)(1.473, 1.562)(1.292, 1.744)(1.141, 1.937)(1.023,
2.134)(0.941, 2.325)(0.909, 2.500)};%
\draw[gen] plot coordinates{(0.913, 2.500)(0.943, 2.665)(1.018,
2.847)(1.128, 3.034)(1.268, 3.220)(1.435, 3.396)(1.627,
3.557)(1.839, 3.699)(2.069, 3.817)(2.312, 3.909)(2.566,
3.971)(2.826, 4.003)(3.087, 4.003)(3.347, 3.971)(3.601,
3.909)(3.844, 3.817)(4.074, 3.699)(4.286, 3.557)(4.478,
3.396)(4.645, 3.220)(4.785, 3.034)(4.895, 2.847)(4.970,
2.665)(5.000, 2.500)(4.970, 2.335)(4.895, 2.153)(4.785,
1.966)(4.645, 1.780)(4.478, 1.604)(4.286, 1.443)(4.074,
1.301)(3.844, 1.183)(3.601, 1.091)(3.347, 1.029)(3.087,
0.997)(2.826, 0.997)(2.566, 1.029)(2.312, 1.091)(2.069,
1.183)(1.839, 1.301)(1.627, 1.443)(1.435, 1.604)(1.268,
1.780)(1.128, 1.966)(1.018, 2.153)(0.943, 2.335)(0.913, 2.500)};%
\draw[gen] plot coordinates{(0.917, 2.500)(0.944, 2.657)(1.015,
2.829)(1.117, 3.008)(1.247, 3.186)(1.402, 3.357)(1.581,
3.515)(1.779, 3.656)(1.993, 3.776)(2.222, 3.873)(2.462,
3.944)(2.708, 3.987)(2.958, 4.001)(3.208, 3.987)(3.455,
3.944)(3.694, 3.873)(3.923, 3.776)(4.138, 3.656)(4.336,
3.515)(4.514, 3.357)(4.670, 3.186)(4.800, 3.008)(4.902,
2.829)(4.972, 2.657)(5.000, 2.500)(4.972, 2.343)(4.902,
2.171)(4.800, 1.992)(4.670, 1.814)(4.514, 1.643)(4.336,
1.485)(4.138, 1.344)(3.923, 1.224)(3.694, 1.127)(3.455,
1.056)(3.208, 1.013)(2.958, 0.999)(2.708, 1.013)(2.462,
1.056)(2.222, 1.127)(1.993, 1.224)(1.779, 1.344)(1.581,
1.485)(1.402, 1.643)(1.247, 1.814)(1.117, 1.992)(1.015,
2.171)(0.944, 2.343)(0.917, 2.500)}; %
\draw[gen] plot coordinates{(0.920, 2.500)(0.946, 2.649)(1.011,
2.813)(1.107, 2.985)(1.228, 3.155)(1.373, 3.320)(1.539,
3.475)(1.724, 3.614)(1.926, 3.736)(2.141, 3.836)(2.367,
3.914)(2.601, 3.966)(2.840, 3.993)(3.080, 3.993)(3.319,
3.966)(3.553, 3.914)(3.779, 3.836)(3.994, 3.736)(4.196,
3.614)(4.381, 3.475)(4.547, 3.320)(4.692, 3.155)(4.813,
2.985)(4.909, 2.813)(4.974, 2.649)(5.000, 2.500)(4.974,
2.351)(4.909, 2.187)(4.813, 2.015)(4.692, 1.845)(4.547,
1.680)(4.381, 1.525)(4.196, 1.386)(3.994, 1.264)(3.779,
1.164)(3.553, 1.086)(3.319, 1.034)(3.080, 1.007)(2.840,
1.007)(2.601, 1.034)(2.367, 1.086)(2.141, 1.164)(1.926,
1.264)(1.724, 1.386)(1.539, 1.525)(1.373, 1.680)(1.228,
1.845)(1.107, 2.015)(1.011, 2.187)(0.946, 2.351)(0.920, 2.500)};%
\definecolor{gen}{rgb}{0,0,1}%
\draw[color=gen,very thick]
plot coordinates{(5.000, 2.500)(4.945, 2.343)(4.871,
2.199)(4.784, 2.065)(4.689, 1.941)(4.586, 1.825)(4.477,
1.719)(4.362, 1.621)(4.241, 1.532)(4.116, 1.452)(3.987,
1.381)(3.854, 1.319)(3.717, 1.266)(3.577, 1.222)(3.435,
1.187)(3.289, 1.163)(3.142, 1.148)(2.992, 1.143)(2.840,
1.149)(2.687, 1.166)(2.533, 1.194)(2.377, 1.235)(2.221,
1.289)(2.064, 1.356)(1.907, 1.438)(1.751, 1.538)(1.596,
1.656)(1.441, 1.797)(1.290, 1.967)(1.141, 2.179)(1.000, 2.500)};%
\draw[color=gen,very thick]
plot coordinates{(5.000, 2.500)(4.945, 2.657)(4.871,
2.801)(4.784, 2.935)(4.689, 3.059)(4.586, 3.175)(4.477,
3.281)(4.362, 3.379)(4.241, 3.468)(4.116, 3.548)(3.987,
3.619)(3.854, 3.681)(3.717, 3.734)(3.577, 3.778)(3.435,
3.813)(3.289, 3.837)(3.142, 3.852)(2.992, 3.857)(2.840,
3.851)(2.687, 3.834)(2.533, 3.806)(2.377, 3.765)(2.221,
3.711)(2.064, 3.644)(1.907, 3.562)(1.751, 3.462)(1.596,
3.344)(1.441, 3.203)(1.290, 3.033)(1.141, 2.821)(1.000, 2.500)};%
\draw[-latex] (-0.750 , 2.500) -- (5.750 , 2.500) node[right] {};
\foreach \x/\xtext in { 0.000, 0.500, ..., 5.000  }%
\draw[shift={(\x,2.500)}] (0,0.025) -- (0,-0.025);%
\foreach \x/\xtext in { 0.000, 0.100, ..., 5.000  }
\draw[shift={(\x,2.500)}] (0,0.0125) -- (0,-0.0125);%
\draw (0 , 2.55) -- (0 , 2.450) node[below] {\FONTSIZE{$-1.25$}}(1 ,
2.55) -- (1 , 2.450) node[below] {\FONTSIZE{$-1$}}(2 , 2.55) -- (2 ,
2.450) node[below] {\FONTSIZE{$-0.75$}}(3 , 2.55) -- (3 , 2.450)
node[below] {\FONTSIZE{$-0.5$}}(4 , 2.55) -- (4 , 2.450) node[below]
{\FONTSIZE{$-0.25$}};%
\draw[-latex] (5.000 , -0.20) -- (5.000 , 5.20) node[above] {};
\foreach \y/\ytext in { 0.000, 0.500, ...,5.000  }%
\draw[shift={(5.000, \y)}] (0.025, 0) -- (-0.025, 0);
\foreach \y/\ytext in { 0.000, 0.100, ..., 5.000  }
\draw[shift={(5.000, \y)}] (0.0125, 0) -- (-0.0125, 0);%
\draw (5.050 , 0) -- (4.950 , 0) node[left]
{\FONTSIZE{$-0.5$}}(5.050 , 1) -- (4.950 , 1) node[left]
{\FONTSIZE{$-0.3$}}(5.050 , 4) -- (4.950 , 4) node[left]
{\FONTSIZE{$0.3$}}(5.050 , 5) -- (4.950 , 5) node[left]
{\FONTSIZE{$0.5$}};
\end{tikzpicture}
\end{center}\vspace*{-.5cm}
\caption{\emph{Distributions of the zeros of $(z+1)\cdots(z+d)-d!=0$
for $d=3,\dots,50$. The zeros approach, as $d$ increases, to the
limiting curve $|z^{-z}(z+1)^{1+z}|=1$ (the blue innermost curve).}}
\label{fig-zeros}
\end{figure}

In terms of the zeros $\lambda_j$'s of the equation $(z+1)
\cdots(z+d)-d!$, we have
\begin{align}\label{d-simplex-mean}
\begin{split}
    \prod_{1\le j<n}
    \left(1-\frac{d!}{(dj+1)\cdots(dj+d)}\right)
    &= \frac{\prod_{1\le j< n}\left(dj\prod_{1\le \ell<d}
    (dj-\lambda_\ell)\right)}
    {\prod_{1\le j< n}\left((dj+1)\cdots(dj+d)\right)}\\
    &= \frac1n \prod_{1\le \ell<d}
    \frac{\Gamma\left(n-\frac{\lambda_\ell}d\right)
    \Gamma\left(1+\frac\ell d\right)}
    {\Gamma\left(n+\frac \ell d\right)
    \Gamma\left(1-\frac{\lambda_\ell}d\right)}\\
    &=: \phi(n).
\end{split}
\end{align}
The zeros $\lambda_j$'s are distributed very regularly as showed
in Figure~\ref{fig-zeros}.

Now we apply the integral representation for the $n$-th finite
difference (called Rice's integrals; see \cite{FS95}) and obtain
\[
    \mu_n
    = -\frac1{2\pi i} \int_{\frac12-i\infty}^{\frac12+i\infty}
    \frac{\Gamma(n+1)\Gamma(-s)}{\Gamma(n+1-s)}\, \phi(s)\dd s.
\]
Note that $\phi(s)$ is well defined and has a simple pole at $s=0$.
The integrand then has a double pole at $s=0$; standard calculations
(moving the line of integration to the left and summing the residue
of the pole encountered) then lead to
\[
    \mu_n
    = \frac1{dH_d}\left(H_n + \sum_{1\le \ell <d}
    \left(\psi\left(-\frac{\lambda_\ell}d\right)-\psi\left(\frac\ell d
    \right)\right) \right) + O\left(n^{-\varepsilon}\right),
\]
where the $O$-term can be made more explicit if needed. Here $H_n =
\sum_{1\le j\le n} 1/j$ denotes the harmonic numbers, and $\psi(z)$
denotes the derivative of $\log\Gamma(z)$. Note that to get this
expression, we used the identity
\[
    \frac{(z+1)\cdots(z+d)-d!}{z}
    = \sum_{1\le j\le d} \frac{d!\Gamma(z+d-j+1)}
    {(d-j+1)!\Gamma(z+1)}.
\]

\paragraph{The probability generating function.} Let $P_n(y) :=
\mathbb{E}[y^{Y_n}]$. Then $P_0(y)=1$ and for $n\ge1$
\[
    P_n(y)
    = y\sum_{0\le k<n} \pi_{n,k} P_k(y).
\]
The same procedure used above leads to
\begin{align*}
    P_n(y)
    &= \sum_{0\le k\le n} \binom{n}{k} (-1)^k
    \prod_{0\le j<k} \left(1-\frac{d!y}
    {(dj+1)\cdots(dj+d)}\right)\\
    &= 1+(y-1)\sum_{1\le k\le n} \binom{n}{k} (-1)^{k-1}
    \prod_{1\le j<k} \left(1-\frac{d!y}
    {(dj+1)\cdots(dj+d)}\right)\qquad(n\ge0).
\end{align*}
Let now $|y-1|$ be close to zero and
\[
    (z+1)\cdots(z+d) - d!y
    = \prod_{1\le \ell \le d} \left(z-\lambda_\ell(y)\right).
\]
Note that the $\lambda_\ell$'s are analytic functions of $y$. Let
$\lambda_d(y)$ denote the zero with $\lambda_d(1)=0$. Then we have
\[
    P_n(y)
    =1- \frac{y-1}{2\pi i} \int_{1-\ve-i\infty}^{1-\ve+i\infty}
    \frac{\Gamma(n+1)\Gamma(-s)}{\Gamma(n+1-s)}\, \phi(s,y)\dd s,
\]
where
\[
    \phi(s,y)
    = \frac{\Gamma\left(s-\frac{\lambda_d(y)}d\right)}
    {\Gamma(s+1)\Gamma\left(1-\frac{\lambda_d(y)}d\right)}
    \prod_{1\le \ell<d}
    \frac{\Gamma\left(s-\frac{\lambda_\ell(y)}d\right)
    \Gamma\left(1+\frac\ell d\right)}
    {\Gamma\left(s+\frac{\ell}d\right)
    \Gamma\left(1-\frac{\lambda_\ell(y)}d\right)}.
\]
Note that for $y\not=1$, $\phi(0,y)=1-y$. When $y\sim 1$, the
dominant zero is $\lambda_d(y)$, and we then deduce that
\[
   P_n(y)
   = Q(y) n^{\lambda_d(y)/d} +O(|1-y|n^{-\ve}),
\]
where
\[
    Q(y)
    := \frac{d(y-1)}
    {\lambda_d(y)\Gamma(1+\frac{\lambda_d(y)}d)}
    \prod_{1\le \ell <d}
    \frac{\Gamma\left(\frac{\lambda_d(y)-\lambda_\ell(y)}
    {d}\right)\Gamma\left(1+\frac\ell d\right)}{\Gamma
    \left(\frac{\lambda_d(y)+\ell}{d}\right)
    \Gamma\left(1-\frac{\lambda_\ell(y)}{d}\right)}.
\]
By writing $(z+1)\cdots(z+d)-d! y =0$ as
\[
    \left(1+z\right)\cdots\left(1+\frac{z}d\right)-1
    = y-1,
\]
and by Lagrange's inversion formula, we obtain
\[
    \lambda_d(y)
    = \frac{y-1}{H_d} -\frac{H_d^2-H_d^{(2)}}{2H_d^3}(y-1)^2
    +O\left(|y-1|^3\right).
\]
From this we then get $Q(1)=1+O(|y-1|)$ and
\[
    \lambda_d(e^\eta)
    = \frac{\eta}{H_d}+\frac{H_d^{(2)}}{2H_d^3}\eta^2 -
    \frac{2H_dH_d^{(3)}-3(H_d^{(2)})^2}{6H_d^5}\eta^3
    +O(|\eta|^4),
\]
for small $|\tau|$. This is a typical situation of the quasi-power
framework (see \cite{FS09,Hwang98}), and we deduce (\ref{mu-Yn}),
(\ref{var-Yn}) and the Berry-Esseen bound (\ref{Yn-BEB}). The
expression for $c_2$ is obtained by an ad-hoc calculation based on
computing the second moment (the expression obtained by the
quasi-power framework being less explicit).

When $d=2$, a direct calculation leads to the identity
\[
    \mathbb{V}[Y_n]
    = \frac{5}{27}H_n +\frac{2\pi^2}{27}+\frac{H_n^{(2)}}{9}
    -\frac{26}{27} -
    \frac29\sum_{j\ge1}\left(\frac{2j-1}{j^2\binom{n+j}{n}}
    -\frac{2j}{(j+\frac12)^2\binom{n+j+\frac12}{n}}\right),
\]
for $n\ge1$, which is also an asymptotic expansion.
This is to be contrasted with $\mathbb{E}[Y_n]=(H_n+2)/3$.

\subsection{Chain records of random samples from hypercubes.}
In this case, we have, denoting still by $Y_n$ the number of chain
records in iud random samples from $[0,1]^d$,
\[
   Y_n \stackrel{d}{=} 1+Y_{I_n} \qquad(n\ge1),
\]
with $Y_0=0$ and
\[
    \mathbb{P}(I_n=k)
    = \binom{n-1}{k} \int_0^1 t^k (1-t)^{n-1-k}
    \frac{(-\log t)^{d-1}}{(d-1)!}\dd t.
\]
Let $P_n(y) := \mathbb{E}[y^{Y_n}]$. Then the Poisson generating
function $\tilde{\mathscr{P}}(z,y) := e^{-z} \sum_{n\ge0} P_n(y)
z^n/n!$ satisfies
\[
    \tilde{\mathscr{P}}(z,y) +
    \frac{\partial}{\partial z} \tilde{\mathscr{P}}(z,y)
    = y\int_0^1 \tilde{\mathscr{P}}(tz)
    \frac{(-\log t)^{d-1}}{(d-1)!} \dd t,
\]
with $\tilde{\mathscr{P}}(0,y) = 1$. We then deduce that
\[
    P_n(y)
    = 1+\sum_{1\le k\le n} \binom{n}{k}
    (-1)^k \prod_{1\le j\le k} \left(1-\frac{y}{j^d}\right).
\]
Consequently, by Rice's integral representation \cite{FS95},
\[
    P_n(y)
    = \frac1{2\pi i}\int_{1-\ve-i\infty}^{1+\ve+i\infty}
    \frac{\Gamma(n+1)\Gamma(-s)}{\Gamma(n+1-s)\Gamma(s+1)^d}
    \prod_{1\le \ell\le d}
    \frac{\Gamma(s+1-y^{1/d} e^{2\ell \pi i/d})}
    {\Gamma(1-y^{1/d} e^{2\ell\pi i/d})} \dd s.
\]
If $|y-1|$ is close to zero, we deduce that
\[
    P_n(y)
    = \frac{n^{y^{1/d}-1}}{\Gamma(y^{1/d})^{1/d}}
    \prod_{1\le\ell<d} \frac{\Gamma(y^{1/d}(1-e^{2\ell\pi i}))}
    {\Gamma(1-y^{1/d}e^{2\ell\pi i/d})}\left(1+O(n^{-\ve})\right).
\]
A very similar analysis as above then leads to a Berry-Esseen bound
for $Y_n$ as follows.

\begin{thm} The number of chain records $Y_n$ for iud random samples
from the hypercube $[0,1]^d$ satisfies \label{thm-clt-cr}
\begin{align*}
    \sup_{x\in\mathbb{R}} \left|\mathbb{P}
    \left(\frac{Y_n-\mu_h\log n}{\sigma_h \sqrt{\log n}}<x\right)
    -\Phi(x)\right| = O\left((\log n)^{-1/2}\right),
\end{align*}
where $\mu_h =\sigma_h := 1/d$. The mean and the variance are
asymptotic to
\begin{align*}
    \mathbb{E}[Y_n]
    &= \frac1d\log n+\gamma
    + \frac1d\sum_{1\le \ell <d}
    \psi(1-e^{2\ell \pi i/d}) +O(n^{-\ve}),\\
    \mathbb{V}[Y_n]
    &= \frac1{d^2} \log n + \frac\gamma d
    -\frac{\pi^2}{6d} \\
    &\quad + \frac1{d^2}\sum_{1\le \ell <d} \left(
    \psi\left(1-e^{2\ell\pi i/d}\right)
    + \left(1-2e^{2\ell\pi i/d}\right)
    \psi'\left(1-e^{2\ell\pi i/d}\right)\right)+O(n^{-\ve}).
\end{align*}
\end{thm}
The asymptotic normality (without rate) was already established in
\cite{Gnedin07}.

In the special case when $d=2$, more explicit expressions are
available
\[
    \mathbb{E}[Y_n]
    = \frac{H_n+1}2,\quad
    \mathbb{V}[Y_n] = \frac{H_n+H_n^{(2)}-2}4,
\]
for $n\ge1$.

\section{Dominating records in the $d$-dimensional simplex}

We consider the mean and the variance of the number of dominating
records in this section.

Let $Z_n$ denote the number of dominating records of $n$ iud points
$\mb{p}_1,\ldots , \mb{p}_n$ in the $d$-dimensional simplex $S_d$.

\begin{thm} The mean and the variance of the number of dominating
records for iud random samples from the $d$-dimensional simplex are
given by
\begin{align}
    \mathbb{E}[Z_n]
    &= \sum_{1\le k\le n}\frac{(d!)^k
    \Gamma(k)^d}{\Gamma(dk+1) },\label{mean-dom-rec} \\
    \mathbb{V}[Z_n]
    &= 2\sum_{2\le k\le n}\frac{(d!)^{k}
    \Gamma(k)^d}{\Gamma(dk+1)}\,H_{k-1}^{(d)}
    +\sum_{1\le k\le n}\frac{(d!)^k
    \Gamma(k)^d}{\Gamma(dk+1) } - \left(
    \sum_{1\le k\le n}\frac{(d!)^k
    \Gamma(k)^d}{\Gamma(dk+1) }\right)^2,\label{var-dom-rec}
\end{align}
respectively. The corresponding expressions for iud random samples
from hypercubes are given by $H_n^{(d)}$ and $H_n^{(d)}-H_n^{(2d)}$,
respectively.
\end{thm}
\pf
\begin{align*}
    \mathbb{E}[Z_n]
    &=\sum_{1\le k\le n}\mathbb{P}\left( \mb{p}_{k}
    \text{ is a dominating record}\right)\\
    &=\sum_{1\le k\le n}(d!)^{k}\int_{S_d}
    \left( \prod_{1\le i\le d} x_i\right) ^{k-1}\dd \mb{x}\\
    &=\sum_{1\le k\le n}\frac{(d!)^{k}}{k}
    \prod_{1\le j<d}\frac{\Gamma(k)\Gamma(jk+1)}{\Gamma((j+1)k+1)}.
\end{align*}
Thus, we obtain (\ref{mean-dom-rec}). For large $n$ and bounded $d$,
the partial sum converges to the series
\[
    \mathbb{E}[Z_n]
    \to \sum_{k\ge1}\frac{(d!) ^{k}
    \Gamma(k)^d}{\Gamma(dk+1)},
\]
at an exponential rate. For large $d$, the right-hand side is 
asymptotic to
\[
    \mathbb{E}[Z_n]
    =1+O\left(\frac{(d!)^2}{(2d)!}\right)
    =1+O\left( 4^{-d}\sqrt d\right),
\]
by Stirling's formula.

Similarly, for the second moment, we have
\begin{align*}
    \mathbb{E}[Z_n^2]-\mathbb{E}[Z_n]
    &=2\sum_{2\le k\le n}\sum_{1\le j<k}\mathbb{P}
    \left( \mb{p}_{j}\text{ and }\mb{p}_{k}
    \text{ are both dominating records}\right) \\
    &=2\sum_{2\le k\le n}\sum_{1\le j<k}
    (d!) ^{k}\int_{S_d}
    \int_{\mb{y}\prec \mb{x}}\left( \prod y_i \right)^{j-1}
    \left( \prod x_i\right) ^{k-j-1}\dd \mb{y}\dd \mb{x}\\
    &=2\sum_{2\le k\le n}(d!)^{k}\int_{S_d}\left(
    \sum_{1\le j<k}\int_{\mb{y}\prec \mb{x}}\left( \prod
    (y_i/x_i)\right)^{j-1}\dd \mb{y}
    \right)\left(\prod x_i\right) ^{k-2}\dd \mb{x} \\
    &=2\sum_{2\le k\le n}(d!)^{k}H_{k-1}^{(d)}
    \int_{S_d}\left( \prod x_i\right)^{k-1}\dd \mb{x} \\
    &=2\sum_{2\le k\le n}\frac{(d!)^{k}
    \Gamma(k)^d}{\Gamma(dk+1)}\,H_{k-1}^{(d)},
\end{align*}
and we obtain (\ref{var-dom-rec}).

For large $n$, the right-hand side of (\ref{var-dom-rec}) converges to
\[
    2\sum_{k\ge 2}\frac{(d!)^{k}
    \Gamma(k)^d}{\Gamma(dk+1)}\,H_{k-1}^{(d)}
    +\sum_{k\ge1}\frac{(d!)^k
    \Gamma(k)^d}{\Gamma(dk+1) } - \left(
    \sum_{k\ge1}\frac{(d!)^k
    \Gamma(k)^d}{\Gamma(dk+1) }\right)^2
\]
at an exponential rate, which, for large $d$, is asymptotic to
$3\sqrt{\pi d}\, 4^{-d}$. This explains the curves corresponding to
$Z_n$ in Figure~\ref{fig-dom-rec}.

The proof for the dominating records in hypercubes is similar and
omitted. \qed


\begin{thebibliography}{99}

\bibitem{BDHT05} Z.-D. Bai, L. Devroye, H.-K. Hwang and H.-T. Tsai,
Maxima in hypercubes, \emph{Random Structures Algorithms}
\textbf{27} (2005), 290--309.

\bibitem{BHLT01} Z.-D. Bai, H.-K. Hwang, W.-Q. Liang and T.-H. Tsai,
Limit theorems for the number of maxima in random samples from
planar regions, \emph{Electron. J. Probab.} \textbf{6} (2001), no.
3, 41 pp. (electronic).

\bibitem{BDJ99} J. Baik, P. Deift and K. Johansson, On the
distribution of the length of the longest increasing subsequence
of random permutations, \emph{J. Amer. Math. Soc.} \textbf{12}
(1999), 1119--1178.

\bibitem{CHT09} W.-M. Chen, H.-K. Hwang and T.-H. Tsai, Simple,
efficient maxima-finding algorithms for multidimensional samples,
preprint, 2009.

\bibitem{CFH07} H.-H. Chern, M. Fuchs and H.-K. Hwang, Phase
changes in random point quadtrees, \emph{ACM Trans. Algorithms}
\textbf{3} (2007), no.\ 2, Art.\ 12, 51 pp.

\bibitem{CQ97} S. N. Chiu and M. P. Quine, Central limit theory
for the number of seeds in a growth model in $\bold R^d$ with
inhomogeneous Poisson arrivals, \emph{Ann. Appl. Probab.} \textbf{7}
(1997),  802--814.

\bibitem{DHBS09} A. Darrasse, H.-K. Hwang, O. Bodini and
M. Soria, The connectivity-profile of random increasing k-trees,
preprint (2009); available at arXiv:0910.3639v1.

\bibitem{DZ95} J.-D. Deuschel and O. Zeitouni, Limiting curves for
iid records, \emph{Ann. Probab.} \textbf{23} (1995), 852--878.

\bibitem{FGD95} P. Flajolet, X. Gourdon and P. Dumas, Mellin
transforms and asymptotics: harmonic sums, \emph{Theoret. Comput.
Sci.} \textbf{144} (1995), 3--58.

\bibitem{FLLS95} P. Flajolet, G. Labelle, L. Laforest and B. Salvy,
Hypergeometrics and the cost structure of quadtrees, \emph{Random
Structures Algorithms} \textbf{7} (1995), 117--144.

\bibitem{FS95} P. Flajolet and R. Sedgewick, Mellin transforms and
asymptotics: finite differences and Rice's integrals, \emph{Theoret.
Comput. Sci.} \textbf{144} (1995), 101--124.

\bibitem{FS09} P. Flajolet and R. Sedgewick, \emph{Analytic
Combinatorics}, Cambridge University Press, Cambridge, 2009.

\bibitem{Gnedin98} A. Gnedin, Records from a multivariate normal
sample, \emph{Statist. Prob. Lett.} \textbf{39} (1998) 11--15.

\bibitem{Gnedin07} A. Gnedin, The chain records, \emph{Electron.
J. Probab.} \textbf{12} (2007), no. 26, 767--786.

\bibitem{Goldie89} C. M. Goldie and S. Resnick, Records in a
partially ordered set, \emph{Ann. Probab} \textbf{17} (1989),
678--699.

\bibitem{Goldie95} C. M. Goldie and S. Resnick, Many multivariate
records, \emph{Stoch. Process. Appl.} \textbf{59} (1995), 185--216.

\bibitem{HH02} E. Hashorva and J. H\"usler, On asymptotics of
multivariate integrals with applications to records, \emph{Stoch.
Models} \textbf{18} (2002), 41--69.

\bibitem{HH05} E. Hashorva and J. H\"usler, Multiple maxima in
multivariate samples, \emph{Stoch. Prob. Letters} \textbf{75}
(2005), 11--17.

\bibitem{Hwang98} H.-K. Hwang, On convergence rates in the central
limit theorems for combinatorial structures, \emph{European J.
Combin.} \textbf{19}  (1998), 329--343.

\bibitem{Kaluszka95} M. Ka\l uszka, Estimates of some probabilities
in multidimensional convex records, \emph{Appl. Math. (Warsaw)}
\textbf{23} (1995), 1--11.

\bibitem{VK77} A. M. Vershik and S. V. Kerov, Asymptotic behavior of
the Plancherel measure of the symmetric group and the limit form of
Young tableaux, \emph{Soviet Math. Dokl.} \textbf{233} (1977),
527--531.

\end{thebibliography}
\end{document}